\documentclass[11pt]{article}
\usepackage[english]{babel}
\usepackage{amsthm}
\usepackage{amsmath}
\usepackage{amssymb}
\usepackage{amsmath,amscd}
\usepackage{tikz-cd}
\usepackage{booktabs}
\usepackage{amsthm}

    \theoremstyle{plain}

\theoremstyle{empty}

\newtheorem{theorem}{Theorem}
\newtheorem{lemma}[theorem]{Lemma}
\newtheorem{proposition}[theorem]{Proposition}
\theoremstyle{proposition*}
\newtheorem*{proposition*}{Proposition}

\theoremstyle{theorem*}

\theoremstyle{theorem*}
\newtheorem*{corollary*}{Corollary}
 
\newtheorem*{theorem*}{Theorem}
\theoremstyle{lemma*}
\newtheorem*{lemma*}{Lemma}
\theoremstyle{definition}

\theoremstyle{definition*}
\newtheorem*{definition*}{Definition}

\theoremstyle{remark}

\begin{document}
\date{{\bf   J Vittis}}
\title{{\bf The D(2)-property for some metacyclic groups}}
\author{}
\maketitle
\setcounter{page}{1}

\begin{abstract}
\noindent 
 
\noindent Keywords  :   Metacyclic group ;   
D(2)-problem. 

\noindent {\it  Mathematics Subject Classification  (AMS  2010):}\,

\noindent {\it  Primary  57M20 :
Secondary 18G50 }
\medskip

We study problems relating to the D(2)-Problem for metacyclic groups of type $G(p,p-1)$ where $p$ is an odd prime.

Specifically we build on Nadim's thesis \cite{Jamil}, which showed that the $\mathbb{Z}[G(5,4)]$-module $\mathbb{Z}$ admits a diagonal resolution and a minimal representative for the third syzygy $\Omega_3(\mathbb{Z})$ is $R(2)\oplus[y-1)$. Motivated by this result, we show that the $\mathbb{Z}[G(p,p-1)]$-module $R(2)\oplus[y-1)$ is both full and straight for any odd prime $p$. Given Johnson's work on the D(2)-Problem \cite{D2}, this leads to the conclusion that $G(5,4)$ satisfies the D(2)-property, as well as providing a sufficient condition for the D(2)-property to hold for $G(p,p-1)$, namely the condition that $R(2)\oplus[y-1)$ is a minimal representative for $\Omega_3(\mathbb{Z})$ over $\mathbb{Z}[G(p,p-1)]$, which we refer to as the condition M(p).

Following this result, we prove a theorem which simplifies the calculations required to show that the condition M(p) holds. Finally, we carry out these calculations in the case where $p=7$ and prove that the condition M(7) holds, which is sufficient to show that $G(7,6)$ satisfies the D(2)-property.
\end{abstract}
\begin{section}{Introduction}

The motivation for this paper comes primarily from Wall's D(2)-Problem which was first formulated in \cite{Wall}.\\

\textbf{The D(2)-Problem:} Let $X$ be a finite connected cell complex with geometric dimension $3$ and with universal cover $\tilde{X}$ such that:
\[
H_{3}(\tilde{X};\mathbb{Z})=0 \text{ and } H^{3}(X;\mathcal{B})=0,
\]
for all coefficient systems $\mathcal{B}$ on $X$. Is $X$ homotopy equivalent to a finite complex of dimension $2$?\\

The D(2)-Problem is intrinsically connected to a second problem in topology, which is known as the two-dimensional realization problem, or the R(2)-Problem.

Let $\mathcal{G}=<x_1,\dots ,x_g \mid W_1, \dots ,W_r>$ be a presentation for a group $G$, and let $K_{\mathcal{G}}$ be the presentation complex of $\mathcal{G}$. Let $\tilde{K{_\mathcal{G}}}$ be the universal cover of $K_{\mathcal{G}}$, also known as the Cayley complex of $\mathcal{G}$. The cellular chain complex of $\tilde{K_{\mathcal{G}}}$ gives rise to:
\[
C_*(\mathcal{G})=(0\rightarrow \pi_2(K_{\mathcal{G}})\rightarrow C_2(\tilde{K_{\mathcal{G}}})\xrightarrow{\partial_2} C_1(\tilde{K_{\mathcal{G}}}) \xrightarrow{\partial_1} C_0(\tilde{K_{\mathcal{G}}}) \xrightarrow{\partial_0} \mathbb{Z} \rightarrow 0),
\]
an exact sequence of right $\mathbb{Z}[G]$-modules.\\

Since each $C_n(\tilde{K_{\mathcal{G}}})$ is a free $\mathbb{Z}[G]$-module, this construction suggests that it might be informative to consider \textit{algebraic 2-complexes} over $G$, which are defined to be exact sequences of $\mathbb{Z}[G]$-modules of the form
\[
0 \rightarrow J \rightarrow F_2 \rightarrow F_1 \rightarrow F_0 \rightarrow \mathbb{Z} \rightarrow 0,
\]
where each $F_n$ is a free $\mathbb{Z}[G]$-module. By the \textit{third syzygy} of $\mathbb{Z}$ over $\mathbb{Z}[G]$, denoted by $\Omega_3(\mathbb{Z})$, we mean the stable module $[J]$, this is well defined as a result of Schanuel's lemma. While considering algebraic 2-complexes, an obvious question arises: whether each algebraic 2-complex over $G$ can be written as $C_*(\mathcal{G})$ where $\mathcal{G}$ is some presentation for $G$.\\

\textbf{The R(2)-Problem:} Let $G$ be a finitely presented group. Is every algebraic 2-complex over $G$
\[
0 \rightarrow J \rightarrow F_2 \rightarrow F_1 \rightarrow F_0 \rightarrow \mathbb{Z} \rightarrow 0,
\]
\textit{geometrically realizable}; that is, homotopy equivalent to an algebraic 2-complex of the form $C_*(\mathcal{G})$, where $\mathcal{G}$ is some presentation for $G$? \\

Johnson showed in \cite{D2} that for finite groups $G$, the R(2)-Problem is equivalent to the D(2)-Problem, that is, if the R(2)-Problem holds true for a finite group $G$, then the D(2)-Problem holds true for all cell complexes $X$ satisfying $\pi_1(X)=G$ and vice versa. This result has since been extended further by Johnson, before reaching its current form, due to Mannan \cite{Mannan}: the R(2)-Problem and the D(2)-Problem are equivalent for all finitely presented groups $G$. If the D(2)-Problem holds true for a group $G$, we say that $G$ satisfies the \textit{D(2)-property}. In this paper, we focus on problems relating to the D(2)-Problem for metacyclic groups of type
\[
G(p,p-1)=<x,y \mid x^p=1, y^{p-1}=1, yx=x^my >,
\]
where $p$ is an odd prime and $m$ is chosen so that the group isomorphism $\theta \in Aut(C_p)$ given by $\theta(x)=x^m$ satisfies $ord(\theta)=p-1$. In Nadim's thesis \cite{Jamil}, some work has already been completed relating to the existence of a diagonal resolution for $\mathbb{Z}$ over $\mathbb{Z}[G(5,4)]$, we use this work as motivation to study the more general case of $G(p,p-1)$.

We begin by stating a sufficient condition for the D(2)-property to hold for a finite group $G$, which follows easily from results in \cite{D2}.

\noindent \textbf{Theorem \ref{the sufficient condition}.}
If a finite group $G$ satisfies properties 1,2 and 3 below
\begin{enumerate}
\item $G$ admits a balanced presentation;
\item $\Omega_3(\mathbb{Z})$ is straight;
\item the minimal module $J$ in $\Omega_3(\mathbb{Z})$ is full;
\end{enumerate}
then $G$ satisfies the D(2)-property.\\

This theorem motivates the problems which are discussed in the remainder of the paper.\\

It is already known \cite{Wamsley} that condition 1 in the above sufficient condition is satisfied in the case where $G=G(p,p-1)$ and $p$ is an odd prime. In \cite{Jamil}, Nadim showed that over $\Lambda=\mathbb{Z}[G(5,4)]$, $\Omega_3(\mathbb{Z})$ has minimal representative $R(2)\oplus [y-1)$, motivated by this result, we aim to show that over $\Lambda=\mathbb{Z}[G(p,p-1)]$, $R(2)\oplus[y-1)$ is full and the stable module $[R(2)\oplus[y-1)]$ is straight. In order to show that these conditions hold, it is useful as a prerequisite to study the ring and $\Lambda$-module
\[
\mathcal{T}_{p-1}(\mathbb{Z},p)=\{(a_{i,j})_{1\leq i,j \leq p-1} \in M_{p-1}(\mathbb{Z}) \mid a_{i,j} \in p \mathbb{Z} \text{ if } i>j \}.
\]
The study of $\mathcal{T}_{p-1}(\mathbb{Z},p)$ is carried out in section \ref{sectiontp-1}, and begins by outlining results from \cite{JJ}, namely a group presentation 
\[
\lambda:G(p,p-1) \rightarrow \mathcal{T}_{p-1}(\mathbb{Z},p),
\]
which extends to a ring surjection
\[
\tilde{\lambda}:\mathbb{Z}[G(p,p-1)]\rightarrow \mathcal{T}_{p-1}.
\]
These results are then used to endow $\mathcal{T}_{p-1}(\mathbb{Z},p)$ with a right $\Lambda$-module structure. As a right $\Lambda$-module, $\mathcal{T}_{p-1}(\mathbb{Z},p)$ is a direct sums of its rows, with this in mind, we denote by $R(i)$ the $i^{th}$ row of $\mathcal{T}_{p-1}(\mathbb{Z},p)$, and so, as right $\Lambda$-modules,
\[
\mathcal{T}_{p-1}(\mathbb{Z},p) \cong \bigoplus_{i=1}^{p-1}R(i).
\]
To conclude the section, we provide a full description of the rings $Hom_{\Lambda}(R(i),R(j))$ and $Hom_{\mathcal{D}er}(R(i),R(j))$ for $1 \leq i,j \leq p-1$.\\

In section \ref{sectionR(2)+[y-1)}, we continue to work over $\Lambda=\mathbb{Z}[G(p,p-1)]$, where $p$ is an odd prime. The section is made up of proofs of two of our main results, mentioned earlier:\\

\noindent \textbf{Theorem \ref{R(2)+[y-1) is straight theorem}.}
$[R(2)\oplus [y-1)]$ is straight over $\Lambda = \mathbb{Z}[G(p,p-1)]$ for any odd prime $p$.\\

\noindent \textbf{Theorem \ref{R(2)+[y-1) is full theorem}.}
$R(2)\oplus [y-1)$ is full over $\Lambda=\mathbb{Z}[G(p,p-1)]$ for any odd prime $p$.\\

In section \ref{G(7,6) chapter)}, we begin by defining the condition \textbf{M(p)} on $\mathbb{Z}[G(p,p-1)]$ as follows:\\

\textbf{M(p):} The third syzygy of $\mathbb{Z}$ over $\Lambda=\mathbb{Z}[G(p,p-1)]$, $\Omega_3(\mathbb{Z})$, is the stable module $[R(2)\oplus[y-1)]$.\\

Given \textbf{Theorem \ref{the sufficient condition}}, \textbf{Theorem \ref{R(2)+[y-1) is straight theorem}}, and \textbf{Theorem \ref{R(2)+[y-1) is full theorem}}, this immediately leads to the following theorem:\\

\noindent \textbf{Theorem \ref{If M(p) holds}.}
Let $\Lambda=\mathbb{Z}[G(p,p-1)]$, if $\Lambda$ satisfies \textbf{M(p)}, then $G(p,p-1)$ satisfies the D(2)-property.\\

As previously noted, it has already been shown \cite{Jamil} that the condition \textbf{M(5)} is satisfied, leading to another one of our main theorems:\\

\noindent \textbf{Theorem \ref{G(5,4) D(2) property}.}
$G(5,4)$ satisfies the D(2)-property\\

The remainder of section \ref{G(7,6) chapter)} is dedicated to refining techniques used in \cite{Jamil} and using these refinements to show that the condition \textbf{M(7)} holds, which leads to our conclusion and final theorems:\\

\noindent \textbf{Theorem \ref{surjection}.}
Over $\Lambda=\mathbb{Z}[G(7,6)]$, $\Omega_3(\mathbb{Z})=[R(2)\oplus [y-1)]$ i.e. the condition $M(7)$ holds.\\

\noindent \textbf{Theorem \ref{syzygy}.}
The D(2)-property holds for $G=G(7,6)$.

\end{section}
\vspace{10mm}
\begin{section}{$\mathcal{T}_{p-1}(\mathbb{Z},p)$}\label{sectiontp-1}
The following is a straightforward consequence of (\cite{D2}, Theorem III):

\begin{theorem}\label{the sufficient condition}
If a finite group $G$ satisfies properties $1,2$ and $3$ below
\begin{enumerate}
\item $G$ admits a balanced presentation;
\item $\Omega_3(\mathbb{Z})$ is straight;
\item the minimal module $J$ in $\Omega_3(\mathbb{Z})$ is full;
\end{enumerate}
then $G$ satisfies the $D(2)$-property.
\end{theorem} 
Let $\Lambda=\mathbb{Z}[G(p,p-1)]$ where $p$ is an odd prime. It is already known \cite{Wamsley} that $G(p,p-1)$ admits a balanced presentation for any odd prime $p$. To utilise the above theorem, it would clearly be useful to describe a minimal module $J \in \Omega_3(\mathbb{Z})$ over $\Lambda$. Currently, such a description only exists in the case where $p=5$. In this case, the minimal module is shown in \cite{Jamil} to be $R(2)\oplus[y-1)$, where $[y-1)$ is the right $\Lambda$-module generated by $[y-1)$, and $R(2)$ is defined as follows:

Let 
\[
\mathcal{T}_{p-1}(\mathbb{Z},p) = \{A=(a_{i,j})_{1,\leq i,j \leq p-1} \in M_{p-1}(\mathbb{Z}) \mid a_{i,j} \in p\mathbb{Z} \text{ if } i>j \}.
\]
For brevity, we will denote $\mathcal{T}_{p-1}(\mathbb{Z},p)$ by $\mathcal{T}_{p-1}$. In \cite{JJ}, a surjective ring homomorphism $\tilde{\lambda}:\Lambda \twoheadrightarrow \mathcal{T}_{p-1}(\mathbb{Z},p)$ is given, which allows us to endow $\mathcal{T}_{p-1}$ with a right $\Lambda$-module structure. Clearly, $\mathcal{T}_{p-1}$ is isomorphic as a right $\Lambda$-module to a direct sum of its rows, which we denote by $R(1),R(2), \dots ,R(p-1)$ from top to bottom, that is:
\[
\mathcal{T}_{p-1} \cong \bigoplus_{i=1}^{p-1}R(i).
\]
We are now motivated to show that over $\Lambda = \mathbb{Z}[G(p,p-1)]$, the stable module $[R(2)\oplus [y-1)]$ is straight and the module $R(2)\oplus [y-1)$ is full. Before proving these properties, we must first discuss the rings $Hom_{\Lambda}(R(i),R(j))$ and $Hom_{\mathcal{D}er}(R(i),R(j))$.
\begin{subsection}{$Hom(R(i),R(j))$}\label{sectionHom(R(i),R(j))}
In this section, we will find explicit descriptions of the rings $Hom_{\Lambda}(R(i),R(j))$, and then state the description of $Hom_{\mathcal{D}er}(R(i),R(j))$ from \cite{JJ}. We begin by defining a mapping $f:\mathbb{Z}^2\rightarrow \{0,p-1\}$ by
\[
    f(i,j)= 
\begin{cases}
    p-1,& \text{if } i>j,\\
    0,              & \text{otherwise.}
\end{cases}
\]
If $\epsilon(i,j)$ is the $(p-1)\times (p-1)$ matrix described by $\epsilon(i,j)_{r,s}=\delta_{i,r}\delta_{j,s}$, where $\delta$ here is the Kronecker delta, then $\mathcal{T}_{p-1}$ has $\mathbb{Z}$-basis given by 
\[
\{t(i,j)=\epsilon(i,j)(1+f(i,j)) \mid 1\leq i,j \leq p-1 \}.
\]
We now define $p-1$ vectors in $M_{1\times (p-1)}(\mathbb{Z})$, each with a right $\Lambda$-action given via $\mathcal{T}_{p-1}$ in the obvious way. 
\begin{itemize}
\item $a_1 =\left(\begin{smallmatrix}1&0&\dots & 0 & 0\end{smallmatrix}\right);$
\item $a_2 =\left(\begin{smallmatrix}0&1&\dots & 0 & 0\end{smallmatrix}\right);$\\
$\vdots$
\item $a_{p-1} =\left(\begin{smallmatrix}0&0&\dots & 0 & 1\end{smallmatrix}\right).$

\end{itemize}
We can think of $R(i)$ as having a $\mathbb{Z}$-basis \[\{a_j(1+f(i,j))=a_it(i,j) \mid 1 \leq j \leq p-1 \}.\] We have shown:
\begin{proposition*}
$R(i)$ is generated over $\Lambda$ by $a_i$.
\end{proposition*}
Therefore, any $\Lambda$-homomorphism $\varphi:R(i)\rightarrow R(j)$ is defined completely by $\varphi(a_i)$. Assume that 
\[
\varphi(a_i)=\sum_{n=1}^{p-1}x_na_n \in R(j).
\]
Where $x_n \in \mathbb{Z}$ for each $n$. Note that we have not placed any further restrictions on the values of $x_n$, and so $\varphi(a_i)$ need not be in $R(j)$ as things stand, this is done in order to simplify the following calculation, and the discrepancy is dealt with shortly. Now, $a_i\sum_{m\neq i} t(m,m) =0$, and so
\[
0=\varphi(a_i\sum_{m \neq i}t(m,m))=\sum_{n=1}^{p-1}x_na_n\sum_{m\neq i}t(m,m)=\sum_{n\neq i }x_na_n,
\]
and so $x_n=0$ for $n \neq i$. Concluding,
\[\varphi(a_i)=x_i a_i. \]
Finally,
\[
\varphi(a_i t(i,k))=\varphi(a_i)t(i,k) = x_ia_it(i,k).
\]
We have shown:
\begin{proposition*}
Let $\varphi \in Hom_{\Lambda}(R(i),R(j))$, then there exists an integer $x_i$ such that
\[ \varphi(\alpha)=\alpha \cdot (x_iI_{p-1}) \; \text{ for each }\alpha \in R(i)\]
\end{proposition*}
We can therefore think of elements of $Hom_{\Lambda}(R(i),R(j))$ as right multiplication by some $nI_{p-1}\in \mathcal{T}_{p-1}$ for some integer $n$. To deal with the discrepancy mentioned above, we must ensure that right multiplication of elements of $R(i)$ by $nI_{p-1}$ gives an element of $R(j)$. This can be ensured by placing a condition on $n$, we clearly have two cases:
\begin{itemize}
\item If $i \geq j$, $Im(nI_{p-1}) \subset R(j)$ for each $n \in \mathbb{Z}$;
\item If $i<j$, $Im(nI_{p-1}) \subset R(j)$ for each $n \in p\mathbb{Z}$.
\end{itemize}
We have shown:
\
\begin{proposition*}
\[
    Hom_{\Lambda}(R(i),R(j))= 
\begin{cases}
    \text{right multiplication by } nI_{p-1},\;n \in \mathbb{Z},& \text{if } i \geq j,\\
     \text{right multiplication by } pnI_{p-1},\;n \in \mathbb{Z}& \text{if } i<j.
\end{cases}
\]
\end{proposition*}

It is clear that $\Lambda$-homomorphisms $\varphi:R(i) \rightarrow R(j)$ and $\psi:R(j)\rightarrow R(k)$ compose in the obvious manner, i.e. if $\varphi = nI_{p-1}$ and $\psi = mI_{p-1}$, then $\psi \circ \varphi = mnI_{p-1}:	R(i)\rightarrow R(k)$. Using this description, it is clear that $End_{\Lambda}(\mathcal{T}_{p-1}(\mathbb{Z},p) )\cong \mathcal{T}_{p-1}(\mathbb{Z},p)$, as expected.

We conclude with a result from \cite{JJ} which describes the rings $Hom_{\mathcal{D}er}(R(i),R(j))$.

\begin{proposition}\label{homder(R(i))}
\[
Hom_{\mathcal{D}er}(R(i),R(j))=
\begin{cases}
\mathbb{Z}/p\mathbb{Z}  &\text{ if } i=j,\\
0  & \text{ if } i \neq j.
\end{cases}
\]
\end{proposition}

\end{subsection}
\end{section}
\begin{section}{$R(2)\oplus [y-1)$}\label{sectionR(2)+[y-1)}
In this section, we will prove that over $\Lambda=\mathbb{Z}[G(p,p-1)]$ where $p$ is an odd prime, the stable module $[R(2)\oplus[y-1)]$ is straight and the module $R(2)\oplus [y-1)$ is full.
\begin{subsection}{Straightness of $R(2)\oplus [y-1)$}\label{sectionstraightnessofR(2)+[y-1)}
\begin{subsubsection}{An exact sequence}
By the Swan-Jacobinski Theorem \cite{C&R2}, in order to show that $[R(2)\oplus [y-1)]$ is straight over $\Lambda = \mathbb{Z}[G(p,p-1)]$, it is sufficient to show that if $S$ is a $\Lambda$-module satisfying
\[
\Lambda \oplus S \cong R(2)\oplus[y-1)\oplus \Lambda,
\]
then $S \cong R(2)\oplus [y-1)$. To do this, we begin by showing that $S$ can be expressed as an extension of $I(C_{p-1})$ by $R(2)\oplus \bigoplus _{i \neq 1}R(i)$ where $I(C_{p-1})$ will be defined shortly. For brevity, we denote $\bigoplus _{i \neq k}R(i)$ by $R(\hat{k})$. Assume that $S$ is a $\Lambda$-module such that 
\[
\Lambda \oplus S \cong R(2)\oplus[y-1)\oplus \Lambda. 
\]
Let $I(C_{p-1})$ be the kernel of the augmentation mapping $\epsilon:\mathbb{Z}[C_{p-1}]\rightarrow \mathbb{Z}$ when considered as a $\Lambda$-homomorphism. Consider the short exact sequence
\[
0\rightarrow R(\hat{1}) \rightarrow [y-1) \rightarrow I(C_{p-1})\rightarrow 0,
\]
from \cite{JJ}. We can now easily construct a second exact sequence
\begin{equation}\label{E sequence}
\mathcal{E}=(0\rightarrow R(2)\oplus R(\hat{1})\rightarrow R(2)\oplus [y-1) \rightarrow I(C_{p-1})\rightarrow 0),
\end{equation}
in the obvious manner. It is also shown in \cite{JJ} that another exact sequence 
\begin{equation}\label{quasiaug}
0 \rightarrow \bigoplus_{j=1}^{p-1}R(j) \rightarrow \Lambda \rightarrow \mathbb{Z}[C_{p-1}] \rightarrow 0,
\end{equation}
exists. By taking the direct sum of this exact sequence with $\mathcal{E}$, we can construct the following exact sequence:
\[
0\rightarrow R(2)\oplus R(\hat{1}) \oplus\bigoplus_{j=1}^{p-1}R(j) \rightarrow R(2)\oplus [y-1)\oplus \Lambda \rightarrow I(C_{p-1})\oplus \mathbb{Z}[C_{p-1}]\rightarrow 0.
\]
By assumption, $S \oplus \Lambda \cong R(2)\oplus [y-1) \oplus \Lambda$, and so we also have a short exact sequence
\begin{equation}\label{to man}
0\rightarrow R(2)\oplus R(\hat{1}) \oplus\bigoplus_{j=1}^{p-1}R(j) \rightarrow S\oplus \Lambda \xrightarrow{p} I(C_{p-1})\oplus \mathbb{Z}[C_{p-1}]\rightarrow 0.
\end{equation}
We will now manipulate (\ref{to man}) to extract an exact sequence $\Phi$ of the following form:
\begin{equation}\label{Phi}
\Phi=(0\rightarrow R(2)\oplus R(\hat{1}) \rightarrow S \rightarrow I(C_{p-1})\rightarrow 0).
\end{equation}
To complete the manipulation we will require some propositions.

\begin{proposition*}
$Hom_{\Lambda}(\mathcal{T}_{p-1},\mathbb{Z}[C_{p-1}]\oplus I(C_{p-1}))=0$
\begin{proof}
As a result of work in \cite{JJ}, as $\Lambda$-modules, $\mathcal{T}_{p-1}(\mathbb{Z},p)\cong [x-1)$. Therefore, for any $t \in \mathcal{T}_{p-1}$, $t\cdot (1+x+\dots +x^{p-1})=0$. Now, as $x$ acts trivially on $\mathbb{Z}[C_{p-1}]$ and $\mathbb{Z}[C_{p-1}]$ is a $\Lambda$-lattice, $Hom_{\Lambda}(\mathcal{T}_{p-1},\mathbb{Z}[C_{p-1}])=0$. By noting that $I(C_{p-1})$ is a $\Lambda$-submodule of $\mathbb{Z}[C_{p-1}]$, the result follows immediately.
\end{proof}
\end{proposition*}

\begin{proposition*}
For the surjection $p:S\oplus \Lambda \rightarrow I(C_{p-1})\oplus \mathbb{Z}[C_{p-1}]$ defined by the short exact sequence \text{(\ref{to man})}, $p(S)\cap p(\Lambda) = 0$.
\begin{proof}
For a $\Lambda$-lattice $M$, we set $M_{\mathbb{Q}}=M \otimes_{\mathbb{Z}}\mathbb{Q}$, which we think of as a $\Lambda_{\mathbb{Q}}$-module. By Wedderburn's Theorem, the isomorphism $S_{\mathbb{Q}} \cong R(2)_{\mathbb{Q}} \oplus [y-1)_{\mathbb{Q}}$ follows from the isomorphism $S\oplus \Lambda \cong R(2) \oplus [y-1)\oplus \Lambda$. By the Wedderburn-Maschke Theorem, we can use the exact sequence (\ref{E sequence}) to form a split exact sequence 
\[
0 \rightarrow (R(2)\oplus R(\hat{1}) )_{\mathbb{Q}} \rightarrow S_{\mathbb{Q}} \rightarrow I(C_{p-1})_{\mathbb{Q}} \rightarrow 0.
\]
Now, $Hom_{\Lambda_{\mathbb{Q}}}((\mathcal{T}_{p-1})_ {\mathbb{Q}},(I(C_{p-1})\oplus\mathbb{Z}[C_{p-1}])_{\mathbb{Q}})=0$, therefore
\[
Hom_{\Lambda_\mathbb{Q}}(S_{\mathbb{Q}},(I(C_{p-1})\oplus \mathbb{Z}[C_{p-1}])_{\mathbb{Q}})\cong Hom_{\Lambda_\mathbb{Q}}(I(C_{p-1})_{\mathbb{Q}},I(C_{p-1})_{\mathbb{Q}}\oplus\mathbb{Z}[C_{p-1}]_{\mathbb{Q}}).\]
We deduce that the image of $\Lambda$-homomorphisms $f:S \rightarrow I(C_{p-1})\oplus \mathbb{Z}[C_{p-1}]$ can have rank of at most $p-2$. By utilising the split exact sequence
\[
0\rightarrow \mathcal{T}_{p-1_{ \mathbb{Q}}} \rightarrow \Lambda_{\mathbb{Q}} \rightarrow \mathbb{Z}[C_{p-1}]_{\mathbb{Q}} \rightarrow 0,
\]
in a similar manner, we see that the maximal rank of the image of a $\Lambda$-homomorphism $g:\Lambda \rightarrow I(C_{p-1})\oplus \mathbb{Z}[C_{p-1}]$ is $p-1$. As $p$ is surjective and $rk_{\mathbb{Z}}(I(C_{p-1})\oplus \mathbb{Z}[C_{p-1}])=(p-2)+(p-1)$ we can now deduce the result.
\end{proof}
\end{proposition*}
Let $J_1=S\cap ker(p)$ and $J_2=\Lambda \cap ker(p)$. We can now think of (\ref{to man}) as a diagonal short exact sequence of the following form.
\[
0 \rightarrow J_1\oplus J_2 \rightarrow S \oplus \Lambda \xrightarrow{p} p(S)\oplus p(\Lambda) \rightarrow 0.
\]
In \cite{Rosen}, it is shown that $\mathcal{T}_{p-1}$ is hereditary, and straightforward calculations using Milnor squares show that $K_0(\mathcal{T}_{p-1})$ is generated by $\{R(j)\}_{1 \leq j \leq p-1}$. Therefore, the modules $J_1$ and $J_2$ must be direct sums of $R(i)$ modules, namely $J_1=\bigoplus_{s=1}^mR(a_s)$ and $J_2 = \bigoplus_{t=1}^n R(b_t)$ for some $a_s,b_t$ such that $m+n=2(p-1)$. Recall the exact sequence
\setcounter{equation}{1}
\begin{equation}
0\rightarrow \bigoplus_{i=1}^{p-1}R(i)\rightarrow \Lambda \rightarrow \mathbb{Z}[C_{p-1}]\rightarrow 0,
\end{equation}
\setcounter{equation}{4}
Consider the following commutative diagram with exact rows:

\begin{center} 
\begin{tikzcd}
 0 \arrow[r] & \bigoplus_{i=1}^{p-1}R(i) \arrow{r}  & \Lambda \arrow{r}    \arrow{d}{Id} & \mathbb{Z}[C_{p-1}]\arrow[r]  & 0 \\  
    0 \arrow[r] & \bigoplus_{t=1}^{n}R(b_t) \arrow{r}  &  \Lambda\arrow{r}  & p(\Lambda) \arrow{r}  & 0.
\end{tikzcd}
\end{center}
Now, as $Hom_{\Lambda}(\mathcal{T}_{p-1},p(\Lambda))=0$ and $Hom_{\Lambda}(\mathcal{T}_{p-1},\mathbb{Z}[C_{p-1}])=0$, $Id:\Lambda \rightarrow \Lambda$ induces and restricts to isomorphisms \[\mathbb{Z}[C_{p-1}]\cong p(\Lambda ),\] and \[\bigoplus_{i=1}^{p-1}R(i)  \cong \bigoplus_{t=1}^{n}R(b_t),\] respectively. Clearly $n=p-1$ and so $m=p-1$. We now wish to show that $p(S) \cong I(C_{p-1})$, to do this we require a proposition.

\begin{proposition*}
Let $J$ be a $\mathbb{Z}[C_{p-1}]$-module such that \[J \oplus \mathbb{Z}[C_{p-1}]\cong I(C_{p-1})\oplus \mathbb{Z}[C_{p-1}],\] as $\mathbb{Z}[C_{p-1}]$-modules. Then $J \cong I(C_{p-1})$.
\begin{proof}
Let $J$ satisfy the above hypothesis. By stabilising the augmentation sequence of $\mathbb{Z}[C_{p-1}]$, we construct the following exact sequence:
\[
0 \rightarrow I(C_{p-1})\oplus \mathbb{Z}[C_{p-1}] \rightarrow \mathbb{Z}[C_{p-1}]^{(2)}\rightarrow \mathbb{Z} \rightarrow 0.
\]
By substituting $J \oplus \mathbb{Z}[C_{p-1}]$ for $I(C_{p-1})\oplus \mathbb{Z}[C_{p-1}]$, we construct the exact sequence 
\[
0 \rightarrow J\oplus \mathbb{Z}[C_{p-1}] \xrightarrow{j} \mathbb{Z}[C_{p-1}]^{(2)}\rightarrow \mathbb{Z} \rightarrow 0.
\]
Let $T=\mathbb{Z}[C_{p-1}]^{(2)}/Im(j|_{\mathbb{Z}[C_{p-1}]})$, taking quotients, we can now construct the exact sequence
\[
0 \rightarrow  J \rightarrow T\rightarrow \mathbb{Z} \rightarrow 0.
\]
By the de-stabilization lemma \cite{syzygies}, $T$ is projective and so the exact sequence \[
0 \rightarrow \mathbb{Z}[C_{p-1}]\xrightarrow{j} \mathbb{Z}[C_{p-1}]^{(2)} \rightarrow T \rightarrow 0,
\] splits, therefore $T \oplus \mathbb{Z}[C_{p-1}] \cong \mathbb{Z}[C_{p-1}]^{(2)}$. Now, $\mathbb{Z}[C_{p-1}]$ satisfies the Eichler condition and so by the Swan-Jacobinski theorem, $\mathbb{Z}[C_{p-1}]$ has stably free cancellation. Therefore, $T\cong \mathbb{Z}[C_{p-1}]$ and we have a short exact sequence 
\[
0 \rightarrow J \rightarrow \mathbb{Z}[C_{p-1}]\rightarrow \mathbb{Z} \rightarrow 0.
\]
Up to sign, the augmentation mapping $\epsilon : \mathbb{Z}[C_{p-1}]\rightarrow \mathbb{Z}$ is the only surjective homomorphism $\mathbb{Z}[C_{p-1}] \rightarrow \mathbb{Z}$ and so $J \cong I(C_{p-1})$ as claimed.
\end{proof}
\end{proposition*}
Using this proposition, we see that $p(S) \cong I(C_{p-1})$ as a $\mathbb{Z}[C_{p-1}]$-module. Given that $p(S)\oplus p(\Lambda) \cong I(C_{p-1})\oplus \mathbb{Z}[C_{p-1}]$ as $\Lambda$-modules, $\mathbb{Z}[C_p]$ clearly acts trivially on $p(S)$, and so $p(S)\cong I(C_{p-1})$ as a $\Lambda$-module. We now have an exact sequence
\[
0\rightarrow \bigoplus_{s=1}^{p-1}R(a_s)\rightarrow S \xrightarrow{p|_S} I(C_{p-1})\rightarrow 0.
\]
Now, $J_2 \cong \bigoplus_{i=1}^{p-1} R(i)$ and $ J_1\oplus J_2 \cong \bigoplus_{s=1}^{p-1}R(a_s)\oplus \bigoplus _{i=1}^{p-1}R(i)\cong R(2)\oplus R(\hat{1})\oplus J_2$, therefore, by \textbf{Proposition \ref{homder(R(i))}.} 
\[
R(2)\oplus R(\hat{1}) \cong   \bigoplus_{s=1}^{p-1}R(a_s).
\]

We have shown that we can express $S$ as an extension $\Phi$ of $I(C_{p-1})$ by $R(2)\oplus R(\hat{1})$:
\setcounter{equation}{3}
\begin{equation}\label{Phi}
\Phi=(0\rightarrow R(2)\oplus R(\hat{1}) \rightarrow S \rightarrow I(C_{p-1})\rightarrow 0).
\end{equation}
Which we wish to compare to
\setcounter{equation}{0}
\begin{equation}\label{E}
\mathcal{E}=(0\rightarrow R(2)\oplus R(\hat{1}) \rightarrow R(2)\oplus [y-1) \rightarrow I(C_{p-1})\rightarrow 0),
\end{equation}
\setcounter{equation}{4}
in order to show that $S \cong R(2)\oplus [y-1)$.
\end{subsubsection}
\begin{subsubsection}{$Ext^1_{\Lambda}(I(C_{p-1}),R(2)\oplus R(\hat{1}))$}\label{Ext R(2)+[y-1)}
In order to compare the extensions (\ref{E}) and (\ref{Phi}), We will find a practical description for $Ext^1_{\Lambda}(I(C_{p-1}),R(2)\oplus R(\hat{1}))$. We will use exact sequences in the derived module category \cite{syzygies} to describe $Ext^1_{\Lambda}(I(C_{p-1}),R(2)\oplus R(\hat{1}))$ as a quotient of the additive abelian group $End_{\mathcal{D}er}(R(2)\oplus R(\hat{1}))$. Recall the short exact sequence
\setcounter{equation}{0}
\begin{equation}\label{E}
\mathcal{E}=(0\rightarrow R(2)\oplus R(\hat{1}) \rightarrow R(2)\oplus [y-1) \rightarrow I(C_{p-1})\rightarrow 0),
\end{equation}
Applying the exact sequence in the derived module category to $\mathcal{E}$, we find an exact sequence
\setcounter{equation}{4}
\begin{equation}\label{long}
Hom_{\mathcal{D}er}(I_{p-1},R(2)\oplus R(\hat{1})) \rightarrow Hom_{\mathcal{D}er}(R(2)\oplus [y-1) ,R(2)\oplus R(\hat{1})) 
\end{equation}
\[
\rightarrow End_{\mathcal{D}er}(R(2)\oplus R(\hat{1})) \rightarrow Ext^1_{\Lambda}(I_{p-1},R(2)\oplus R(\hat{1}))
\]
\[
\rightarrow Ext^1_{\Lambda}(R(2)\oplus [y-1), R(2)\oplus R(\hat{1}))
\rightarrow \dots \]
where, for brevity, we denote $I(C_{p-1})$ by $I_{p-1}$. We now find practical descriptions for some of the groups in this exact sequence, beginning with $Hom_{\mathcal{D}er}(I_{p-1},R(2)\oplus R(\hat{1}))$. As a straightforward consequence of the isomomorphism $[x-1) \cong \mathcal{T}_{p-1}$,  $Hom_{\Lambda}(I_{p-1},R(2)\oplus R(\hat{1}))=0$, therefore $Hom_{\mathcal{D}er}(I_{p-1},R(2)\oplus R(\hat{1}))=0$.\\

We now state without proof a proposition which will allow us to simplify the exact sequence (\ref{long}). 
\begin{proposition*}(\cite{JJ}, page 47, 5.8)
\[
    Ext^1_{\Lambda}(I_{p-1},R(k))= 
\begin{cases}
    0,& \text{if } k=1\\
    \mathbb{Z}/p\mathbb{Z},  & \text{otherwise}
\end{cases}
\]
\end{proposition*}
To calculate the second term, we require a lemma.

\begin{lemma}\label{homder(y-1,t_{p-1}}
$Hom_{\mathcal{D}er}([y-1),\mathcal{T}_{p-1})=0$
\begin{proof}
Recall that $[y-1)$ occurs in an exact sequence 
\[
0 \rightarrow R(\hat{1}) \rightarrow [y-1) \rightarrow I_{p-1} \rightarrow 0.
\]
Applying the exact sequence in the derived module category to the above exact sequence gives rise to a second exact sequence,
\[
Hom_{\mathcal{D}er}(I_{p-1},\mathcal{T}_{p-1}) \rightarrow Hom_{\mathcal{D}er}([y-1),\mathcal{T}_{p-1}) \rightarrow Hom_{\mathcal{D}er}(R(\hat{1}),\mathcal{T}_{p-1})
\]
\[
\rightarrow Ext^1_{\Lambda}(I_{p-1},\mathcal{T}_{p-1}) \rightarrow Ext^1_{\Lambda}([y-1),\mathcal{T}_{p-1}) \rightarrow Ext^1_{\Lambda}(R(\hat{1}),\mathcal{T}_{p-1}).
\]
By the above proposition, $Ext^1_{\Lambda}(I_{p-1},\mathcal{T}_{p-1})=(\mathbb{Z}/p\mathbb{Z})^{p-2}$, using this result in conjunction with \textbf{Proposition \ref{homder(R(i))}} the above exact sequence becomes
\[
0 \rightarrow Hom_{\mathcal{D}er}([y-1),\mathcal{T}_{p-1}) \rightarrow (\mathbb{Z}/p\mathbb{Z})^{p-2}  \rightarrow (\mathbb{Z}/p\mathbb{Z})^{p-2}\]\[ \rightarrow Ext^1_{\Lambda}([y-1),\mathcal{T}_{p-1}) \rightarrow \dots
\]
Therefore, if we can show that $Ext^1_{\Lambda}([y-1),\mathcal{T}_{p-1})=0$, our proof of the lemma will be complete. Let $j:\mathbb{Z}[C_p] \rightarrow \Lambda$ be the standard inclusion, the result  $Ext^1_{\Lambda}([y-1),\mathcal{T}_{p-1})=0$ follows easily from isomorphism $j^*([y-1))\cong \mathbb{Z}[C_{p}]^{p-2}$ and the Eckmann-Shapiro lemma. This completes the proof.
\end{proof}
\end{lemma}

Collecting our results, we can rewrite the exact sequence (\ref{long}) as
\[
0 \rightarrow (\mathbb{Z}/p\mathbb{Z})^2 \rightarrow (\mathbb{Z}/p\mathbb{Z})^{p+1} \rightarrow (\mathbb{Z}/p\mathbb{Z})^{p-1} \xrightarrow{p^*} Ext^1_{\Lambda}(R(2)\oplus [y-1),R(2)\oplus R(\hat{1})) \rightarrow \dots
\]
We deduce that $p^*=0$, rewriting (\ref{long}) a final time, we find the exact sequence
\[
0 \rightarrow Hom_{\mathcal{D}er}(R(2)\oplus [y-1) ,R(2)\oplus R(\hat{1})) 
\xrightarrow{i^*} End_{\mathcal{D}er}(R(2)\oplus R(\hat{1}))
\]
\[ \xrightarrow{\delta} Ext^1_{\Lambda}(I_{p-1},R(2)\oplus R(\hat{1})) \rightarrow 0.
\]
We will now calculate $Im(i^*)$, and use this information to give a practical description for  $Ext^1_{\Lambda}(I_{p-1},R(2)\oplus R(\hat{1}))$. When considered as a matrix, elements $f$ in $End_{\Lambda}(R(2)\oplus R(\hat{1}))$ take the form
\[ f = 
\left( 
\begin{array}{@{}c|c@{}}
     f_1:R(2)\oplus R(2) \rightarrow R(2) \oplus R(2) 
      & f_2:\bigoplus_{i \neq 1,2}R(i)\rightarrow R(2)\oplus R(2)\\
   \cmidrule[0.4pt]{1-2}
   f_3:R(2)\oplus R(2) \rightarrow \bigoplus_{i \neq 1,2}R(i)& f_4:\bigoplus_{i \neq 1,2}R(i)\rightarrow \bigoplus_{i \neq 1,2}R(i)
\end{array}
\right).
\]
By \textbf{Proposition \ref{homder(R(i))}},  the element represented by $f$ in the derived module category, $\overline{f}$ takes the form
\[ \overline{f} = 
\left( 
\begin{array}{@{}c|c@{}}
     \overline{f_1}
      & 0\\
   \cmidrule[0.4pt]{1-2}
   0 & \overline{f_4}
\end{array}
\right)=
\left( 
\begin{array}{@{}c|c@{}}
   \begin{array}{@{}cccc@{}}
      a_{1,1}+p\mathbb{Z} & a_{1,2}+p\mathbb{Z} \\
      a_{2,1}+p\mathbb{Z} & a_{2,2}+p\mathbb{Z}
    
   \end{array} 
      & 0\\
   \cmidrule[0.4pt]{1-2}
   0 & \begin{array}{@{}cccc@{}}
      a_{3,3}+p\mathbb{Z} &  &  & \\
       & a_{4,4}+p\mathbb{Z} & & \\
       & &\ddots & \\
       &  &  & a_{p-1,p-1}+p\mathbb{Z}
    
   \end{array}   \\
\end{array}
\right),
\]
Here, elements in the matrix are zero if not otherwise specified. The elements $a_{n,n}+p\mathbb{Z}$ are in $ End_{\mathcal{D}er}(R(n))$ for $3 \leq n \leq p-1$. To find $Im(i^*)$, we take a general element $f^{\prime} \in Hom_{\Lambda}(R(2)\oplus [y-1),R(2)\oplus R(\hat{1}))$ which, when considered as a matrix takes the form
\[ f^{\prime} = 
\left( 
\begin{array}{@{}c|c@{}}
     f_1^{\prime}:R(2) \rightarrow R(2) \oplus R(2) 
      & f_2^{\prime}:[y-1)\rightarrow R(2)\oplus R(2)\\
   \cmidrule[0.4pt]{1-2}
   f_3^{\prime}:R(2) \rightarrow \bigoplus_{i \neq 1,2}R(i)& f_4^{\prime}:[y-1)\rightarrow \bigoplus_{i \neq 1,2}R(i)
\end{array}
\right).
\]
We proved earlier in this section that $Hom_{\mathcal{D}er}([y-1),\mathcal{T}_{p-1})=0$, and so the element represented by $f^{\prime}$ in the derived module category, $\overline{f}^{\prime}$ is given by
\[ \overline{f}^{\prime} = 
\left( 
\begin{array}{@{}c|c@{}}
     \overline{f_1}^{\prime}:R(2) \rightarrow R(2) \oplus R(2) 
      & \overline{0}:[y-1)\rightarrow R(2)\oplus R(2)\\
   \cmidrule[0.4pt]{1-2}
   \overline{0}:R(2) \rightarrow \bigoplus_{i \neq 1,2}R(i)& \overline{0}:[y-1)\rightarrow \bigoplus_{i \neq 1,2}R(i)
\end{array}
\right)= 
\left( 
\begin{array}{@{}c|c@{}}
     \overline{f_1}^{\prime}
      & 0\\
   \cmidrule[0.4pt]{1-2}
   0 & 0
\end{array}
\right).
\]
Now, $i^*(\overline{f^{\prime}})=\overline{f^{\prime} \circ i}=\overline{f^{\prime}} \circ \overline{i}$ and so 
\[
Im(i^*)=\{
\left( 
\begin{array}{@{}c|c@{}}
        \begin{array}{@{}cccc@{}}
      a_{1,1}+p\mathbb{Z} & 0+p\mathbb{Z} \\
      a_{2,1}+p\mathbb{Z} & 0+p\mathbb{Z}
    \end{array}
      & 0\\
   \cmidrule[0.4pt]{1-2}
   0 & 0
\end{array}
\right)\in End_{\mathcal{D}er}(R(2)\oplus R(\hat{1})) \mid a_{1,1},a_{2,1}\in \mathbb{Z}   \}.
\]
We conclude that given an $\alpha$ in $End_{\Lambda}(R(2)\oplus R(\hat{1}))$, which takes the form
\[
\left(
\begin{array}{cccccccc}

a_{1,1}&a_{1,2}&a_{1,3}&a_{1,4}&\dots&a_{1,p-2}&a_{1,p-1}\\
a_{2,1}&a_{2,2}&a_{2,3}&a_{2,4}&\dots&a_{2,p-2}&a_{2,p-1}\\
pa_{3,1}&pa_{3,2}&a_{3,3}&a_{3,4}&\dots&a_{3,p-2}&a_{3,p-1}\\
pa_{4,1}&pa_{4,2}&pa_{4,3}&a_{4,4}&\dots&a_{4,p-2}&a_{4,p-1}\\
\vdots & \vdots & \vdots & \vdots & \ddots & \vdots & \vdots \\
pa_{p-1,1}&pa_{p-1,2}&pa_{p-1,3}&pa_{p-1,4}&\dots&pa_{p-1,p-2}&a_{p-1,p-1}
\end{array}
\right),
\]
the element in $End_{\mathcal{D}er}(R(2)\oplus \bigoplus_{i \neq 1}R(i)) / Im(i^*)$ represented by $\alpha$, which we denote by $[\alpha]$ is given by
\[
[\alpha]=
\left( 
\begin{array}{@{}c|c@{}}
   \begin{array}{@{}cccc@{}}
      0 & a_{1,2}+p\mathbb{Z} \\
      0 & a_{2,2}+p\mathbb{Z}
    
   \end{array} 
      & 0\\
   \cmidrule[0.4pt]{1-2}
   0 & \begin{array}{@{}cccc@{}}
      a_{3,3}+p\mathbb{Z} &  &  & \\
       & a_{4,4}+p\mathbb{Z} & & \\
       & &\ddots & \\
       &  &  & a_{p-1,p-1}+p\mathbb{Z}
    
   \end{array}   \\
\end{array}
\right).
\]
Using our rewritten version of (\ref{long}), we know that \[End_{\mathcal{D}er}(R(2)\oplus R(\hat{1})/Im(i^*)\cong Ext^1_{\Lambda}(I_{p-1},R(2)\oplus R(\hat{1}))\] via the mapping $[\alpha] \mapsto \alpha_*(\mathcal{E})$, we have shown:

\begin{theorem}\label{Ext^1}
Every element of $Ext^1_{\Lambda}(I_{p-1},R(2)\oplus R(\hat{1}))$ can be written uniquely as $\alpha_*(\mathcal{E})$, where $\alpha \in End_{\Lambda}(R(2)\oplus R(\hat{1}))$ is given by a matrix
\[
\alpha=
\left( 
\begin{array}{@{}c|c@{}}
   \begin{array}{@{}cccc@{}}
      0 & a_{1,2}\\
      0 & a_{2,2}
    
   \end{array} 
      & 0\\
   \cmidrule[0.4pt]{1-2}
   0 & \begin{array}{@{}cccc@{}}
      a_{3,3} &  &  & \\
       & a_{4,4} & & \\
       & &\ddots & \\
       &  &  & a_{p-1,p-1}
    
   \end{array}   \\
\end{array}
\right),
\]
where $a_{i,j} \in \mathbb{Z}$ and $0 \leq a_{i,j} \leq p-1$ for each $i,j$. Elements of the matrix are zero if not otherwise specified.
\end{theorem}
Recall $\Phi$, the short exact sequence (\ref{Phi}), we may now deduce that $\Phi=\alpha_*(\mathcal{E})$ for some $\alpha \in End_{\Lambda}(R(2)\oplus R(\hat{1}))$ of the above form. We now prove some lemmas which allow us to place restrictions on the form of $\alpha$.

\begin{lemma*}
Assume that $\alpha_*(\mathcal{E})=\Phi$, then it can not be true that $a_{1,2}=a_{2,2}=0$.
\begin{proof}
Recall that $\Phi$ takes the form\[
\Phi=(0\rightarrow R(2)\oplus \bigoplus_{i \neq 1} R(i) \rightarrow S \rightarrow I(C_{p-1})\rightarrow 0).
\]
Assume that $a_{1,2}=a_{2,2}=0$, then basic properties of the colimit imply that
\[
S \cong R(2) \oplus R(2) \oplus X,
\]
for some $\Lambda$-module $X$. We know that $S \oplus \Lambda \cong R(2)\oplus [y-1) \oplus \Lambda$, but \[
Hom_{\mathcal{D}er}(R(2)\oplus R(2) \oplus X\oplus \Lambda, R(2)) \cong (\mathbb{Z}/p\mathbb{Z})^2 \oplus Hom_{\mathcal{D}er}(X,R(2)),
\]
while
\[
Hom_{\mathcal{D}er}(R(2)\oplus [y-1)\oplus \Lambda, R(2)) \cong \mathbb{Z}/p\mathbb{Z}.
\]
This gives a contradiction, completing the proof.
\end{proof}
\end{lemma*}
Similarly, we can prove the following:
\begin{lemma*}
Assuming that $\alpha_*(\mathcal{E})=\Phi$, it can not be true that $a_{i,i}=0$ for any $i$ satisfying $3 \leq i \leq p-1$.
\end{lemma*}
\end{subsubsection}
\begin{subsubsection}{Realising the matrices as automorphisms}
In this section, we will prove that there must exist an automorphism \[\alpha: R(2)\oplus R(\hat{1})\rightarrow R(2)\oplus R(\hat{1}),\] such that $\Phi = \alpha_*(\mathcal{E})$. This result, used in conjunction with the five lemma leads to the conclusion that $R(2)\oplus [y-1) \cong S$. Let $\alpha$ be a general element of $End_{\Lambda}(R(2)\oplus R(\hat{1}))$, when considered as a matrix, $\alpha$ takes the form
\[
\left(
\begin{array}{cccccccc}

a_{1,1}&a_{1,2}&a_{1,3}&a_{1,4}&\dots&a_{1,p-2}&a_{1,p-1}\\
a_{2,1}&a_{2,2}&a_{2,3}&a_{2,4}&\dots&a_{2,p-2}&a_{2,p-1}\\
pa_{3,1}&pa_{3,2}&a_{3,3}&a_{3,4}&\dots&a_{3,p-2}&a_{3,p-1}\\
pa_{4,1}&pa_{4,2}&pa_{4,3}&a_{4,4}&\dots&a_{4,p-2}&a_{4,p-1}\\
\vdots & \vdots & \vdots & \vdots & \ddots & \vdots & \vdots \\
pa_{p-1,1}&pa_{p-1,2}&pa_{p-1,3}&pa_{p-1,4}&\dots&pa_{p-1,p-2}&a_{p-1,p-1}
\end{array}
\right).
\]Define a map \[k:End_{\Lambda}(R(2)\oplus \bigoplus_{i \neq 1}R(i))\rightarrow (\mathbb{Z}/p\mathbb{Z})^{p-1}, 
\] by
\[
\alpha \mapsto (a_{1,2}+p\mathbb{Z},a_{2,2}+p\mathbb{Z},a_{3,3}+p\mathbb{Z},a_{4,4}+p\mathbb{Z}, \dots ,a_{p-1,p-1}+p\mathbb{Z}).
\]
By \S \ref{Ext R(2)+[y-1)}, $End_{\mathcal{D}er}(R(2)\oplus R(\hat{1}))/Im(i^*) \cong Ext^1_{\Lambda}(I_{p-1},R(2)\oplus R(\hat{1}))$ via the pushout map $[\alpha ]\mapsto \alpha_*(\mathcal{E})$, and so if  $\alpha,\beta \in End_{\Lambda}(R(2)\oplus R(\hat{1}))$, then $[\alpha]=[\beta]$ if and only if $k(\alpha)=k(\beta)$. Therefore $\alpha_*(\mathcal{E})=\beta_*(\mathcal{E})$ if and only if $k(\alpha)=k(\beta)$. We noted in the lemmas at the end of the previous section that $\Phi = \alpha_*(\mathcal{E})$ for some $\alpha$ such that $k(\alpha)$ can have at most one zero entry, and $a_{i,i}+p\mathbb{Z}\neq 0 +p\mathbb{Z}$ for $3 \leq i \leq p-1$, therefore, if we can find an automorphism $\alpha$ for each of these cases, it must be true, by the five lemma, that $S \cong R(2)\oplus [y-1)$.\\

Define the set of units 
\begin{align*}
U_{(a_{1,2},a_{2,2},\dots,a_{p-1,p-1)}}=\{ & \alpha \in Aut_{\Lambda}(R(2)\oplus R(\hat{1})) \mid  k(\alpha)= \\ & \; \; (a_{1,2}+p\mathbb{Z},a_{2,2}+p\mathbb{Z},a_{3,3}+p\mathbb{Z},\dots ,a_{p-1,p-1}+p\mathbb{Z})\}.
\end{align*}
To prove that $S \cong R(2)\oplus [y-1)$, we will show that $U_{(a_{1,2},a_{2,2},\dots,a_{p-1,p-1)}}$ is non-empty whenever $a_{1,2},a_{2,2}$ are not both zero $(mod p)$ and $a_{i,i} \not\equiv 0 (modp)$ for $3 \leq i \leq p-1$. We begin by noting some generating elements:

\[
f_{(n,1,1,\dots ,1)}=
\left( 
\begin{array}{@{}c|c@{}}
   \begin{array}{@{}cccc@{}}
      1 & n\\
      0 & 1
     \end{array} 
      & 0    \\
   \cmidrule[0.4pt]{1-2}
   0
    & \begin{array}{@{}cccc@{}}
     I_{p-3}
   \end{array}   \\
\end{array}
\right)\in U_{(n,1,1,\dots ,1)},
\]
\[
f_{(1,0,1,\dots ,1)}=
\left( 
\begin{array}{@{}c|c@{}}
   \begin{array}{@{}cccc@{}}
      0 & 1\\
      1 & 0
     \end{array} 
      & 0    \\
   \cmidrule[0.4pt]{1-2}
   0
    & \begin{array}{@{}cccc@{}}
     I_{p-3}
   \end{array}   \\
\end{array}
\right)\in U_{(1,0,1,\dots ,1)},
\]
\[
f_{(0,2,1,1,\dots ,1)}=
\left( 
\begin{array}{@{}c|c@{}}
   \begin{array}{@{}cccc@{}}
      \frac{p+1}{2} & 0\\
      -1 & 2
    
   \end{array} 
      & \begin{array}{@{}cccc@{}}
      1 & 0 & \dots  & 0 \\
      0 & 0 & \dots & 0 \\
    
   \end{array}\\
   \cmidrule[0.4pt]{1-2}
   \begin{array}{@{}cccc@{}}
      0 & p  \\
      0 & 0  \\
      \vdots & \vdots \\
      0 & 0
    
   \end{array}
    & \begin{array}{@{}cccc@{}}
      1 &  &  & \\
       & 1 & & \\
       & &\ddots & \\
       &  &  & 1
    
   \end{array}   \\
\end{array}
\right)\in U_{(0,2,1,\dots ,1)}.
\]
Note that 
\[
U_{(a_{1,2},a_{2,2}, \dots ,a_{p-1,p-1})}\cdot f_{(0,2,1,\dots ,1)} \subset U_{(2a_{1,2},2a_{2,2},a_{3,3},\dots a_{p-1,p-1})}.
\]
Therefore, by considering automorphisms of the form
\[f_{(n,1,1,\dots ,1)}\cdot f_{(0,2,1,1,\dots ,1)}^a \text{ and } f_{(1,0,1,\dots ,1)}\cdot f_{(0,2,1,1,\dots ,1)}^a,
\]
one sees easily that
$U_{(a_{1,2},a_{2,2},1,1\dots ,1)}
\text{ is non-empty whenever }a_{1,2},a_{2,2} $ are not both zero.

We now note two more generating elements:
\[
f_{(0,1,2,1,\dots ,1)}=
\left( 
\begin{array}{@{}c|c@{}}
   \begin{array}{@{}cccc@{}}
      \frac{p+1}{2} & 0\\
      0 & 1
    
   \end{array} 
      & \begin{array}{@{}cccc@{}}
      1 & 0 & \dots  & 0 \\
      0 & 0 & \dots & 0 \\
    
   \end{array}\\
   \cmidrule[0.4pt]{1-2}
   \begin{array}{@{}cccc@{}}
      p & 0  \\
      0 & 0  \\
      \vdots & \vdots \\
      0 & 0
    
   \end{array}
    & \begin{array}{@{}cccc@{}}
      2 &  &  & \\
       & 1 & & \\
       & &\ddots & \\
       &  &  & 1
    
   \end{array}   \\
\end{array}
\right)\in U_{(0,1,2,1,\dots ,1)},
\]

\[
f_{(0,1,\frac{p+1}{2},2,1,\dots ,1)}=
\left( 
\begin{array}{@{}c|c@{}}
   \begin{array}{@{}cccc@{}}
      1 & 0\\
      0 & 1
    
   \end{array} 
      & 0\\
   \cmidrule[0.4pt]{1-2}
   0
    & \begin{array}{@{}ccccc@{}}
      \frac{p+1}{2} & 1  &  & \\
      p  & 2 & & &\\
       & &1 &  &\\
       &  &  & \ddots & \\
       &  &  &  & 1
   \end{array}   \\
\end{array}
\right)\in U_{(0,1,\frac{p+1}{2},2,1,\dots ,1)} .
\]
By considering automorphisms of the form $f_{(0,1,2,1,\dots ,1)}^a$, we see that
\[
U_{(0,1,a_{3,3},1,\dots ,1)} \text{ is non-empty whenever } a_{3,3} \not\equiv 0 (modp).
\]
Now, note that $f_{(0,1,2,1,\dots ,1)}\cdot f_{(0,1,\frac{p+1}{2},2,1,\dots 1)} \in U_{(0,1,1,2,1,\dots ,1)} $. If we define the automorphism $f_{(0,1,1,2,1,\dots ,1)} = f_{(0,1,2,1,\dots ,1)}\cdot f_{(0,1,\frac{p+1}{2},2,1,\dots 1)}$, then by considering automorphisms of the form $f_{(0,1,1,2,1,\dots ,1)}^a$ we can see similarly that 
\[
U_{(0,1,1,a_{4,4},\dots ,1)} \text{ is non-empty whenever } a_{4,4} \not\equiv 0 (modp).
\] 
Repeating this process, we see that $U_{(0,1,1,\dots a_{i,i}, \dots , 1)}$ is non-empty whenever $a_{i,i} \not\equiv 0 (modp)$ for $3 \leq i \leq p-1$. Now, note that 
\[
U_{(0,1,a_{3,3},1,\dots ,1)}\cdot U_{(0,1,1,a_{4,4},\dots ,1)} \cdot \dots \cdot U_{(0,1,1,\dots ,a_{p-1,p-1})} \subset U_{(0,1,a_{3,3},a_{4,4},\dots ,a_{p-1,p-1})}.
\]
We deduce that \[U_{(0,1,a_{3,3},a_{4,4},\dots ,a_{p-1,p-1})} \text{ is non-empty whenever } a_{i,i} \not\equiv 0 (modp)\; 3 \leq i \leq p-1.\]
Finally, we note that 
\[f_{(a_{1,2},a_{2,2},1,\dots ,1)}\cdot U_{(0,1,a_{3,3},a_{4,4},\dots ,a_{p-1,p-1})} \subset U_{(a_{1,2},a_{2,2},a_{3,3},a_{4,4},\dots ,a_{p-1,p-1})}.
\]
We have shown: 
\begin{proposition*}
$U_{(a_{1,2},a_{2,2},\dots,a_{p-1,p-1)}}$ is non-empty whenever $a_{1,2},a_{2,2}$ are not both zero $(mod p)$ and $a_{i,i} \not\equiv 0 (modp)$ for $3 \leq i \leq p-1$
\end{proposition*} Therefore, there exists an $\alpha \in Aut_{\Lambda}(R(2)\oplus \bigoplus_{i \neq 1}R(i))$ such that $\Phi = \alpha_*(\mathcal{E})$. In conclusion:
\begin{theorem}\label{R(2)+[y-1) is straight theorem}
$S \cong R(2)\oplus [y-1)$ and so $[R(2)\oplus [y-1)]$ is straight over $\Lambda = \mathbb{Z}[G(p,p-1)]$ for any odd prime $p$.
\end{theorem}
\end{subsubsection}
\end{subsection}
\begin{subsection}{$R(2)\oplus [y-1)$ is full}\label{sectionfullnessofR(2)+[y-1)}
In this section, we will show that $R(2)\oplus [y-1)$ is full over $\Lambda=\mathbb{Z}[G(p,p-1)]$ for any odd prime $p$. We will begin by showing that $R(2)$ and $[y-1)$ are both full $\Lambda$-modules, before building upon these results to show that $R(2)\oplus [y-1)$ is full.\\

To show that $R(2)$ is full, we must find the kernel of the Swan homomorphism \cite{syzygies} $S_{R(2)}:Aut_{\mathcal{D}er}(R(2))\rightarrow \tilde{K}_0(\Lambda)$. Recall that $End_{\Lambda}(R(2))=\{nId_{R(2)}\mid n \in \mathbb{Z} \} \cong \mathbb{Z}$, and  $End_{\mathcal{D}er}(R(2))\cong \mathbb{Z}/p\mathbb{Z}$, so $Aut_{\mathcal{D}er}(R(2))=\{\overline{nId_{R(2)}}\mid 1 \leq n < p \} \cong (\mathbb{Z}/p\mathbb{Z})^*$, the units of $\mathbb{Z}/p\mathbb{Z}$. Recall the exact sequence
\setcounter{equation}{1}
\begin{equation}
0 \rightarrow \bigoplus_{i=1}^{p-1}R(i)\xrightarrow{\iota} \Lambda \rightarrow \mathbb{Z}[C_{p-1}] \rightarrow 0,
\end{equation}
\setcounter{equation}{5}
Using this exact sequence, we can form a second exact sequence
\[
0 \rightarrow R(2) \xrightarrow{\iota|_{R(2)}} \Lambda \rightarrow \Lambda/\iota(R(2)) \rightarrow 0.
\]
In order to find $Ker(S_{R(2)})$, we must answer the question: 'for which $\overline{n} \in (\mathbb{Z}/p\mathbb{Z})^*$ is $\varinjlim(\overline{nId_{R(2)}}, \iota|_{R(2)})$ stably free?' By the Swan-Jacobinski Theorem \cite{C&R2}, a consequence of $\Lambda$ satisfying the Eichler condition is that $\Lambda$ satisfies SFC i.e. each stably free $\Lambda$-module is free, therefore, it is sufficient to find when $\varinjlim(\overline{nId_{R(2)}}, \iota|_{R(2)}) \cong \Lambda$. Using the Five lemma, it is immediately obvious that for $\overline{n}=\overline{1}$, $\overline{n}=\overline{p}-\overline{1}$, $\varinjlim(\overline{nId_{R(2)}},\iota|_{R(2)})\cong \Lambda$. We will show that these are the only two choices for $\overline{n} \in (\mathbb{Z}/p\mathbb{Z})^*$ such that this isomorphism holds.\\

Let the mapping 
\[\widetilde{nId_{R(2)}}:\bigoplus_{i=1}^{p-1}R(i)\rightarrow \bigoplus_{i=1}^{p-1}R(i),\] be defined by
\[
(r_1,r_2,r_3,\dots ,r_{p-1}) \mapsto (r_1,nr_2,r_3,\dots r_{p-1}),
\]
clearly, $\varinjlim(\widetilde{nId_{R(2)}},\iota) \cong \varinjlim(nId_{R(2)},\iota|_{R(2)})$. From work done in \cite{edwards}, the exact sequence \ref{quasiaug} gives rise to a fibre square \cite{milnor} of the following form

\[ \begin{tikzcd}
\Lambda \arrow{r}{i_1} \arrow[swap]{d}{i_2} & End_{\Lambda}(\mathbb{Z}[C_{p-1}]) \arrow{d}{j_1} \\%
End_{\Lambda}(\bigoplus_{i=1}^{p-1}R(i))\arrow{r}{j_2}& End_{\mathcal{D}er}(\bigoplus_{i=1}^{p-1}R(i)),
\end{tikzcd}
\]
where $i_1,j_1$ and $j_2$ are surjective. The maps $i_1,i_2$ are define as follows: take a $\lambda \in \Lambda$, define $g_{\lambda}:\Lambda \rightarrow \Lambda$ by $g_{\lambda}(1)=\lambda$. $Hom_{\Lambda}(\mathcal{T}_{p-1},\mathbb{Z}[C_{p-1}])=0$, and so $g_{\lambda}$ gives rise to a commutative diagram with exact rows
\begin{center} 
\begin{tikzcd}
 0 \arrow[r] & \bigoplus_{i=1}^{p-1}R(i) \arrow{r}\arrow{d}{\tilde{g}}  & \Lambda \arrow{r}    \arrow{d}{g_{\lambda}} & \mathbb{Z}[C_{p-1}]\arrow[r]\arrow{d}{g}  & 0 \\  
    0 \arrow[r] & \bigoplus_{i=1}^{p-1}R(i) \arrow{r}  &  \Lambda\arrow{r}  &  \mathbb{Z}[C_{p-1}] \arrow{r}  & 0.
\end{tikzcd}
\end{center}
In this case, $i_1(\lambda)=g$ and $i_2(\lambda)=\tilde{g}$. In matrix form, $End_{\Lambda}(\bigoplus_{i=1}^{p-1}R(i)) \cong \mathcal{T}_{p-1}(\mathbb{Z},p)$ and $End_{\mathcal{D}er}(\bigoplus_{i=1}^{p-1}R(i))$ is the ring of $(p-1)\times (p-1)$ matrices with elements of $\mathbb{Z}/p\mathbb{Z}$ on the diagonal, and zeroes elsewhere. The map $j_1:End_{\Lambda}(\bigoplus_{i=1}^{p-1}R(i))\rightarrow End_{\mathcal{D}er}(\bigoplus_{i=1}^{p-1}R(i))$, when considered in matrix form is given by
\setlength{\arraycolsep}{1pt}
\[
\left(
\begin{array}{cccccccc}

a_{1,1}&a_{1,2}&\dots&a_{1,p-1}\\
pa_{2,1}&a_{2,2}&\dots&a_{2,p-1}\\
\vdots & \vdots & \ddots & \vdots &  \\
pa_{p-1,1} & pa_{p-1,2} & \dots & a_{p-1,p-1} &  
\end{array}
\right) \mapsto 
\left(
\begin{array}{cccccccc}

a_{1,1}+p\mathbb{Z}& & & &\\
&a_{2,2}+p\mathbb{Z}& & &\\
& &\ddots&&\\
&& & a_{p-1,p-1}+p\mathbb{Z}&  &   
\end{array}
\right).
\]
By (\cite{edwards}, Proposition 4.5.2), $\varinjlim(\widetilde{nId_{R(2)}},\iota) \cong M(End_{\Lambda}(\mathbb{Z}[C_{p-1}]),\mathcal{T}_{p-1},\widetilde{nId_{R(2)}})$. Projective modules of type  $M(End_{\Lambda}(\mathbb{Z}[C_{p-1}]),\mathcal{T}_{p-1},h)$ are in $1-1$ correspondence with the quotient set 
\[
j_2(End_{\Lambda}(\mathbb{Z}[C_{p-1}])^*)\backslash Aut_{\mathcal{D}er}(\bigoplus_{i=1}^{p-1}R(i))/j_1(\mathcal{T}_{p-1}^*).
\]
\cite{syzygies}. To show that $\varinjlim(\widetilde{nId_{R(2)}},\iota) \ncong \Lambda$ when $n \not\equiv \pm 1(modp)$ it is therefore sufficient to show that the class of $\widetilde{Id_{R(2)}}$ is different to that of $\widetilde{nId_{R(2)}}$ in the above quotient set. In order to do this, we begin by simplifying our description of the quotient set.
\setlength{\arraycolsep}{6pt}

\begin{lemma*}
\[j_2(End_{\Lambda}(\mathbb{Z}[C_{p-1}])^*)\backslash Aut_{\mathcal{D}er}(\bigoplus_{i=1}^{p-1}R(i))/j_1(\mathcal{T}_{p-1}^*) \cong Aut_{\mathcal{D}er}(\bigoplus_{i=1}^{p-1}R(i))/j_1(\mathcal{T}_{p-1}^*)\]
\begin{proof}
Take a unit $f \in End_{\Lambda}(\mathbb{Z}[C_{p-1}])$, $f$ clearly lifts to a unit $f_{\lambda} \in End_{\Lambda}(\Lambda)$, which in turn lifts to an $\tilde{f} \in End_{\Lambda}(\bigoplus_{i =1}^{p-1}R(i))$ such that
\begin{center} 
\begin{tikzcd}
 0 \arrow[r] & \bigoplus_{i=1}^{p-1}R(i) \arrow{r}\arrow{d}{\tilde{f}}  & \Lambda \arrow{r}    \arrow{d}{f_{\lambda}} & \mathbb{Z}[C_{p-1}]\arrow[r]\arrow{d}{f}  & 0 \\  
    0 \arrow[r] & \bigoplus_{i=1}^{p-1}R(i) \arrow{r}  &  \Lambda\arrow{r}  &  \mathbb{Z}[C_{p-1}] \arrow{r}  & 0,
\end{tikzcd}
\end{center}
commutes. By the Five lemma, $\tilde{f}$ is an isomorphism, and by our construction of the Milnor square, $j_1(f)=j_2(\tilde{f})$. Therefore, there exists an $\tilde{f} \in \mathcal{T}_{p-1}^*$ such that $j_1(f)=j_2(\tilde{f})$, completing the proof.
\end{proof}
\end{lemma*}
To show that $\varinjlim(\widetilde{nId_{R(2)}},\iota) \ncong \Lambda$ for $n \not\equiv \pm 1 (modp)$, it is now sufficient to show that no unit in $\mathcal{T}_{p-1}^*$ has $(1(modp),n(modp),1(modp),\dots ,1(modp))$ on the diagonal. Basic considerations relating to the determinant show that any such matrix would have determinant $n(modp)$, and so if $n \not\equiv \pm 1 (modp)$, $\varinjlim(\widetilde{nId_{R(2)}},\iota) \ncong \Lambda$. We have shown:
\begin{proposition*}
Over $\Lambda=\mathbb{Z}[G(p,p-1)]$, $Ker(S_{R(2)})=\{\pm Id\}$.
\end{proposition*}
Now, $\{\pm Id\} \subset Im(v^{R(2)})$ and so $Im(v^{R(2)})=Ker(S_{R(2)})$, concluding:
\begin{proposition*}
Over $\Lambda = \mathbb{Z}[G(p,p-1)]$, $R(2)$ is full.
\end{proposition*}
We will now show that $[y-1)$ is full, to do this, we will first consider the problem over $\mathbb{Z}[C_{p-1}]$. Let $[y-1)^{\prime}$ be the right $\mathbb{Z}[C_{p-1}]$-module generated by $(y-1)$.
\begin{proposition}\label{ker(sy-1)}
Over $\mathbb{Z}[C_{p-1}]$, $[y-1)^{\prime}$ is full and \[S_{[y-1)^{\prime}}:Aut_{\mathcal{D}er}([y-1)^{\prime})\rightarrow \tilde{K}_0(\mathbb{Z}[C_{p-1}])\] is the zero map.
\begin{proof}
Consider the augmentation sequence over $\mathbb{Z}[C_{p-1}]$,
\[
0\rightarrow [y-1)^{\prime} \rightarrow \mathbb{Z}[C_{p-1}]\xrightarrow{\epsilon} \mathbb{Z} \rightarrow 0,
\]
where $\epsilon$ is defined by $\epsilon(y)=1$. We know \cite{syzygies} that $Aut_{\mathcal{D}er}([y-1)^{\prime}) \cong Aut_{\mathcal{D}er}(\mathbb{Z})$, therefore $Aut_{\mathcal{D}er}([y-1)^{\prime})=\{\overline{nId}_{[y-1)^{\prime}}\mid gcd(n,p-1)=1; 1 \leq n \leq p-2 \}$. We will find a representative for each element of $Aut_{\mathcal{D}er}([y-1)^{\prime})$ which is also an element of $Aut_{\Lambda}([y-1)^{\prime})$, and our result will follow by the Five lemma.

Fix an $r$ such that $1 \leq r \leq p-2$ and $gcd(r,p-1)=1$. Let
\[
(1+y+\dots + y^{r-1}):[y-1)^{\prime} \rightarrow [y-1)^{\prime}
\]
be the $\mathbb{Z}[C_{p-1}]$-homomorphism given by  multiplication by $(1+y+\dots + y^{r-1})$. 

Firstly, we will show that $(1+y+\dots +y^{r-1}):[y-1)^{\prime}\rightarrow [y-1)^{\prime}$ is an automorphism by showing that $y-1$ lies in its image. By assumption, $r$ is coprime to $p-1$, so there exists an $s$ such that $rs\equiv 1 (mod(p-1))$. Clearly $y^r-1,y^{2r}-y^r, \dots ,y^{sr}-y^{(s-1)r}$ are elements of $Im(1+y+\dots +y^{r-1})$, and their sum is $y-1$. Therefore $(1+y+\dots +y^{r-1})\in Aut_{\Lambda}([y-1)^{\prime})$.

To complete the proof, it remains only to show that \[\overline{(1+y+\dots +y^{r-1})}=\overline{rId}\in End_{\mathcal{D}er}([y-1)^{\prime}).\] To prove this, it is sufficient to show that
\[
-r+(1+y+\dots +y^{r-1})\in [y-1)^{\prime}.
\]
But
\begin{align*}
-r+(1+y+\dots +y^{r-1}) &= -(r-1)+y+y^2+\dots + y^{r-1},\\
&= (y-1)+(y^2-1)+\dots + (y^{r-1}-1) \in [y-1)^{\prime}.
\end{align*}
Therefore, each $\overline{f}\in Aut_{\mathcal{D}er}([y-1)^{\prime})$ lifts to an automorphism over $\mathbb{Z}[C_{p-1}]$ and $S_{[y-1)^{\prime}}=0$, as required.
\end{proof}
\end{proposition}
Given that over $\Lambda$, $End_{\mathcal{D}er}([y-1))=\{\overline{nId}_{[y-1)}\mid 0 \leq n \leq p-2 \}\cong \mathbb{Z}/(p-1)\mathbb{Z}$, this proposition extends naturally to show:

\begin{proposition}\label{ker(sy-1lambda)}
Over $\Lambda$, $[y-1)$ is full and \[S_{[y-1)}:Aut_{\mathcal{D}er}([y-1))\rightarrow \tilde{K}_0(\mathbb{Z}[C_{p-1}])\] is the zero map.
\end{proposition}
We have shown that both $R(2)$ and $[y-1)$ are full $\Lambda$-modules. In order to extend these results to $R(2)\oplus [y-1)$, we require a lemma.
\begin{lemma}
Let $A,B$ be full $\Lambda$-modules. If $A,B$ satisfy the following properties:
\begin{itemize}
\item $Hom_{\mathcal{D}er}(A,B)=0$;
\item $Hom_{\mathcal{D}er}(B,A)=0$;
\item $S_B: Aut_{\mathcal{D}er}(B)\rightarrow \tilde{K}_0(\Lambda)$ is the zero mapping.
\end{itemize}
Then $A \oplus B$ is full.
\begin{proof}
Let $\overline{f}$ be an element of $Ker(S_{A\oplus B})$, then
\[
\overline{f}=
\begin{pmatrix}
\overline{f_1}:A\rightarrow A&\overline{f_2}:B\rightarrow A \\ 
\overline{f_3}:A\rightarrow B&\overline{f_4}:B\rightarrow B
\end{pmatrix} = 
\begin{pmatrix}
\overline{f_1} & \overline{0} \\
\overline{0} & \overline{f_4}
\end{pmatrix}\in Aut_{\mathcal{D}er}(A\oplus B),
\]
So $\overline{f_1}\in Aut_{\mathcal{D}er}(A)$ and $\overline{f_4}\in Aut_{\mathcal{D}er}(B)$. Now, $S_{A\oplus B}(\overline{f}) = S_A(\overline{f_1})+ S_B(\overline{f_4})$ , by our hypothesis, $S_B=0$ and so $S_{A\oplus B}(\overline{f}) = S_A(\overline{f_1})$. We assumed that $\overline{f} \in Ker(S_{A\oplus B})$ and so $S_A(\overline{f_1})=0$, but $A$ is full and so $\overline{f_1} \in Im(v^A)=Ker(S_A)$. Finally, $\overline{f_4} \in Aut_{\mathcal{D}er}(B)$, but $S_B=0$ and $B$ is full, therefore $\overline{f}_4\in Ker(S_B)=Im(v^B)$. finally, $\overline{f}\in Im(v^{A\oplus B})$ because $\overline{f_1} \in Im(v^A)$ and $\overline{f}_4 \in Im(v^B)$, this completes the proof.
\end{proof}
\end{lemma}

By \textbf{Lemma \ref{homder(y-1,t_{p-1}}} and \textbf{Proposition \ref{ker(sy-1lambda)}}, $Hom_{\mathcal{D}er}([y-1),R(2))=0$ and $S_{[y-1)}=0$. To show that $R(2)\oplus [y-1)$ is full, it remains only to prove that $Hom_{\mathcal{D}er}(R(2),[y-1))=0$, we instead prove the following, stronger result:
\begin{proposition*}
$Hom_{\mathcal{D}er}(\mathcal{T}_{p-1},[y-1))=0.$
\begin{proof}
Recall the exact sequence
\[
0\rightarrow \mathcal{T}_{p-1}\rightarrow \Lambda \rightarrow \mathbb{Z}[C_{p-1}] \rightarrow 0,
\]
from \cite{JJ}. By applying the exact sequence in the derived module category \cite{syzygies}, we see that 
\[
Hom_{\mathcal{D}er}(\mathcal{T}_{p-1},[y-1))\cong Ext^1_{\Lambda}(\mathbb{Z}[C_{p-1}],[y-1)).
\]
The result now follows easily from the Eckmann-Shapiro lemma.
\end{proof}
\end{proposition*}
Collecting our results:

\begin{theorem}\label{R(2)+[y-1) is full theorem}
$R(2)\oplus [y-1)$ is full over $\Lambda=\mathbb{Z}[G(p,p-1)]$.
\end{theorem}
\end{subsection}
\end{section}
\begin{section}{The condition M(7)}\label{G(7,6) chapter)}
In this section, we will begin by defining a condition $M(p)$ over $\Lambda=\mathbb{Z}[G(p,p-1)]$ where $p$ is an odd prime, which we will use to give a practical sufficient condition for the $D(2)$-property to hold for $G(p,p-1)$. We will then prove a theorem which significantly shortens the calculations necessary to show that the condition $M(p)$ holds. We will close the section by proving that the condition $M(7)$ holds, which in turn leads to one of our main results, namely that the D(2)-property holds for the group $G(7,6)$.
\begin{subsection}{The condition \textbf{M(p)}}

Let $\Lambda=\mathbb{Z}[G(p,p-1)]$ where $p$ is an odd prime, we define the condition $\textbf{M(p)}$ on $\Lambda$ as follows:

\indent $\textbf{M(p)}:$ The third syzygy of $\mathbb{Z}$ over $\Lambda=\mathbb{Z}[G(p,p-1)]$ is the stable module $[R(2)\oplus [y-1)]$.

Basic considerations relating to the rank of $R(2)\oplus [y-1)$ as a $\mathbb{Z}$-module show that if $\textbf{M(p)}$ is satisfied, then $R(2)\oplus[y-1)$ is in fact a minimal representative for $\Omega_3^{G(p,p-1)}(\mathbb{Z})$. Therefore, if \textbf{M(p)} is satisfied, our results relating to $R(2)\oplus[y-1)$ in section \ref{sectionR(2)+[y-1)} mean that all three conditions in \textbf{Theorem \ref{the sufficient condition}} are satisfied, we have shown:
\begin{theorem}\label{If M(p) holds}
Let $\Lambda=\mathbb{Z}[G(p,p-1)]$, if $\Lambda$ satisfies \textbf{M(p)}, then $G(p,p-1)$ satisfies the D(2)-property.
\end{theorem}
It has already been shown \cite{Jamil} that $\textbf{M(5)}$ is satisfied, this leads to one of our main theorems.
\begin{theorem}\label{G(5,4) D(2) property}
$G(5,4)$ satisfies the D(2)-property
\end{theorem}
The remainder of this paper is devoted to showing that $\textbf{M(7)}$ is satisfied.

\end{subsection}

\begin{subsection}{A theorem relating to the condition $M(p)$}\label{M(p)}
We wish to study $Ext^1_{\Lambda}(\mathbb{Z}[C_{p-1}],R(\hat{n}))$ by expressing its elements as pushouts of elements of the abelian group $Ext^1_{\Lambda}(\mathbb{Z}[C_{p-1}],\bigoplus_{i=1}^{p-1}R(i))$. Recall the exact sequence,
\[\Psi=(
0\rightarrow \bigoplus_{i=1}^{p-1}R(i) \xrightarrow{j}\Lambda \rightarrow \mathbb{Z}[C_{p-1}]\rightarrow 0)
\]
from \cite{JJ}, which we now denote by $\Psi$. Applying the exact sequence in the derived module category to $\Psi$, we find that 
\[
Hom_{\mathcal{D}er}(\bigoplus_{i=1}^{p-1}R(i),R(\hat{n}))\cong Ext^1_{\Lambda}(\mathbb{Z}[C_{p-1}],R(\hat{n})),
\]
via the map 
\begin{align*}\delta_*:Hom_{\mathcal{D}er}(\bigoplus_{i=1}^{p-1}R(i),R(\hat{n}))&\rightarrow Ext^1_{\Lambda}(\mathbb{Z}[C_{p-1}],R(\hat{n})),\\
\overline{f} &\mapsto f_*(\Psi).
\end{align*}

We can explicitly describe the additive abelian groups $Hom_{\Lambda}(\bigoplus_{i=1}^{p-1}R(i),R(\hat{n}))$ and $Hom_{\mathcal{D}er}(\bigoplus_{i=1}^{p-1}R(i),R(\hat{n}))$ in matrix form:\\

Let $f \in Hom_{\Lambda}(\bigoplus_{i=1}^{p-1}R(i),R(\hat{n}))$, we can think of $f$ as a $(p-2)\times (p-1)$ matrix such that:
\[
    (f)_{i,j} \in 
\begin{cases}
   Hom_{\Lambda}(R(j),R(i)),& \text{if } i \leq n-1 \\
    Hom_{\Lambda}(R(j),R(i+1)),              & \text{if }i \geq n.
\end{cases}
\]
This is equivalent to
\[
    (f)_{i,j} \in 
\begin{cases}
   \mathbb{Z},& \text{if } i \leq n-1\text{ and } j \geq i \\
    p\mathbb{Z},& \text{if } i \leq n-1\text{ and } i \geq j+1 \\
    \mathbb{Z},  & \text{if }i \geq n \text{ and } j \geq i+1\\
      p\mathbb{Z},  & \text{if }i \geq n \text{ and } i+1 \geq j+1.
\end{cases}
\]
We can then express $\overline{f}$ as a $(p-2)\times (p-1)$ matrix such that
\[
    (\overline{f})_{i,j} \in 
\begin{cases}
   \mathbb{Z}/p\mathbb{Z},& \text{if } i=j \text{ and } i \leq n-1 \\
    \mathbb{Z}/p\mathbb{Z},              & \text{if } i+1=j \text{ and } i \geq n\\
    \{0\},& \text{otherwise}.
\end{cases}
\]
Note that the standard projection map \[Hom_{\Lambda}(\bigoplus_{i=1}^{p-1}R(i),R(\hat{n})) \rightarrow Hom_{\mathcal{D}er}(\bigoplus_{i=1}^{p-1}R(i),R(\hat{n}))\] takes the obvious form when considered in matrix form.\\

Consider an $f \in Hom_{\Lambda}(\bigoplus_{i=1}^{p-1}R(i),R(\hat{n}))$ in matrix form, we define an additive group homomorphism
\[
k:Hom_{\Lambda}(\bigoplus_{i=1}^{p-1}R(i),R(\hat{n}))\rightarrow (\mathbb{Z},p\mathbb{Z})^{p-2},
\]
by
\[
f \mapsto (f_{1,1}+p\mathbb{Z},f_{2,2}+p\mathbb{Z},\dots ,f_{n-1,n-1}+p\mathbb{Z},f_{n,n+1}+p\mathbb{Z}, \dots ,f_{p-2,p-1}+p\mathbb{Z}).
\]
Note that $k$ descends to an isomorphism in the derived module category i.e. $k(f)=k(f^{\prime})$ if and only if $\overline{f}=\overline{f^{\prime}}$. Therefore $k(f)=k(f^{\prime})$ if and only if $\delta_*(f)=\delta_*(f^{\prime})$. It is now clear that for each $\Phi \in Ext^1_{\Lambda}(\mathbb{Z}[C_{p-1}],R(\hat{n}))$, there exists an $f \in Hom_{\Lambda}(\bigoplus_{i=1}^{p-1}R(i),R(\hat{n}))$ such that $\Phi = \delta_*(f)$, moreover $\delta_*(f)=\delta_*(f^{\prime})$ if and only if $k(f)=k(f^{\prime})$, we therefore classify $Ext^1_{\Lambda}(\mathbb{Z}[C_{p-1}],R(\hat{k}))$ by the \textit{$k$-invariants} of $f$, which we define to be
\[
f_{1,1}+p\mathbb{Z},f_{2,2}+p\mathbb{Z},\dots ,f_{n-1,n-1}+p\mathbb{Z},f_{n,n+1}+p\mathbb{Z}, \dots ,f_{p-2,p-1}+p\mathbb{Z}.
\]
We aim to show that if there exists a short exact sequence
\[
\mathcal{S}=(0\rightarrow R(\hat{n})\rightarrow K \rightarrow \mathbb{Z}[C_{p-1}]\rightarrow 0),
\]
with all non-zero $k$-invariants, then there exists an imbedding $i:R(n)\hookrightarrow \Lambda$ such that $K\cong \Lambda/i(R(n))$.\\

We begin by defining the set $U_{(x_1,x_2,\dots ,x_{n-1},x_{n+1},x_{n+2}, \dots x_{p-1})}$ as follows:
\begin{align*}
\{ u  \in \mathcal{T}_{p-1}(\mathbb{Z},p)^* \mid & u \text{ has } (y_1,y_2,\dots ,y_{n-1},\star ,y_{n+1}, \dots , y_{p-1}) \text{ on the diagonal }\\
& \text{ and } y_i \equiv x_i (modp) \text{ for all } i\},
\end{align*}
where $\star$ can represent any integer.
\begin{lemma*}
For any $(p-2)$-tuple of integers $(c_1,\dots c_{n-1},c_{n+1},\dots ,c_{p-1})$ such that $c_i \not\equiv0(modp)$ for each $i$, the set $U_{(c_1,\dots c_{n-1},c_{n+1},\dots ,c_{p-1})}$ is non-empty.
\begin{proof}
By B\'ezout's lemma, $U_{(c_1,1,1,\dots)} ,U_{(1,c_2,1,\dots , 1)},\dots U_{(1,1,\dots ,1,c_{p-1})}$ are all non-empty. By taking the product of an element from each of these sets, we construct an element of $U_{(c_1,\dots c_{k-1},c_{k+1},\dots ,c_{p-1})}$, completing the proof.
\end{proof}
\end{lemma*}
Fix a $(p-2)$-tuple of integers $(d_1,d_2,\dots ,d_{n-1},d_{n+1},\dots ,d_{p-1})$ such that $d_i \not\equiv 0(modp)$ for each $i$. We will find an $f \in Hom_{\Lambda}(\bigoplus_{i=1}^{p-1}R(i),R(\hat{n}))$ such that \[k(f) = (d_1+p\mathbb{Z},d_2+p\mathbb{Z},\dots ,d_{n-1}+p\mathbb{Z},d_{n+1}+p\mathbb{Z},\dots ,d_{p-1}+p\mathbb{Z}),\] and $\delta_*(f)$ takes the form 
\[
\delta_*(f)=(0 \rightarrow R(\hat{n}) \rightarrow \Lambda/R(n) \rightarrow \mathbb{Z}[C_{p-1}] \rightarrow 0.
\]  By the lemma, $U_{(d_1,d_2,\dots ,d_{n-1},d_{n+1},\dots ,d_{p-1})}$ is non-empty, and so contains an element $u$. Define $h \in Hom_{\Lambda}(\bigoplus_{i=1}^{p-1}R(i),R(\hat{n}))$ to be the mapping 
\begin{align*}h:R(\hat{n})\oplus R(n) &\rightarrow R(\hat{n})\\
(\hat{n},n) & \mapsto \hat{n}.
\end{align*}

If  $i:\bigoplus_{i=1}^{p-1}R(i) \hookrightarrow \Lambda$ is an injection, basic properties of the colimit give an isomorphism $\varinjlim(h,i)\cong \Lambda/i(R(n))$. This leads to the conclusion that $\varinjlim(h\circ u ,j)\cong \Lambda /i(R(n))$. But 
\[
k(h\circ u) = (d_1+p\mathbb{Z},d_2+p\mathbb{Z},\dots ,d_{n-1}+p\mathbb{Z},d_{n+1}+p\mathbb{Z},\dots ,d_{p-1}+p\mathbb{Z}).
\]
Therefore, if we define $f = h \circ u$, then
\[
k(f)=(d_1+p\mathbb{Z},d_2+p\mathbb{Z},\dots ,d_{n-1}+p\mathbb{Z},d_{n+1}+p\mathbb{Z},\dots ,d_{p-1}+p\mathbb{Z}),
\]
and
\[\delta_*(f) = 0\rightarrow R(\hat{k})\rightarrow \Lambda/i(R(n))\rightarrow \mathbb{Z}[C_{p-1}]\rightarrow 0. \]

Concluding:

\begin{theorem}\label{k-invariants}
For every extension
\[
O \rightarrow R(\hat{n})\rightarrow K \rightarrow \mathbb{Z}[C_{p-1}]\rightarrow 0,
\]
with all non-zero $k$-invariants, there exists an imbedding $i:R(n)\rightarrow \Lambda$ such that \[K \cong \Lambda/i(R(n)).\]
\end{theorem}

This theorem is useful in showing that the condition $M(p)$ holds, as if we can show that a surjective homomorphism $\pi: \Lambda \rightarrow R(1)$ exists with kernel $K$, and that $K$ is an extension of $\mathbb{Z}[C_{p-1}]$ by $R(\hat{2})$ with all non-zero $k$-invariants, then $K \cong \Lambda/R(2)$. Given that the Kernel of the augmentation map $\epsilon:\Lambda \rightarrow \mathbb{Z}$ is isomorphic to $R(1)\oplus [y-1)$ \cite{JJ}, the isomorphism $K \cong \Lambda/R(2)$ then leads to the conclusion that $\Omega_3^{G(p,p-1)}(\mathbb{Z})=[R(2)\oplus [y-1)]$. This is the scheme of proof which was used in \cite{Jamil} to show that the condition $M(5)$ holds. We will use this scheme to show that the condition $M(7)$ holds.
\end{subsection}

\begin{subsection}{$\mathcal{T}_6(\mathbb{Z},7)$}\label{T_6}
The remainder of this paper is dedicated to proving that the condition $M(7)$ holds. Fix the presentation \[G(7,6)=<x,y \mid y^6=1,x^7=1,yx=x^3y>,\] for the group $G(7,6)$. Let $\Lambda$ be the integral group ring of $G(7,6)$. In this section, by studying the representation $\lambda:G(7,6)\rightarrow \mathcal{T}_6(\mathbb{Z},7)^*$ described in \cite{JJ}, we will find characteristic equations for each of the rows of $\mathcal{T}_6(\mathbb{Z},7)$. We will then use these characteristic equations to find representations for each of the rows of $\mathcal{T}_6(\mathbb{Z},7)$. For brevity, we abbreviate $\mathcal{T}_6(\mathbb{Z},7)$ to $\mathcal{T}_6$.

One can calculate $\lambda(x^{-1})$ and $\lambda(y^{-1})$ by hand, giving the following result:

\[ \lambda(x^{-1})=  
 \left(
\begin{array}{cccccc}
 -6 & 1 & 0 & 0 & 0 & 0 \\
 -21 & 1 & 1 & 0 & 0 & 0 \\
 -35 & 0 & 1 & 1 & 0 & 0 \\
 -35 & 0 & 0 & 1 & 1 & 0 \\
 -21 & 0 & 0 & 0 & 1 & 1 \\
 -7 & 0 & 0 & 0 & 0 & 1 \\
\end{array}
\right),
\]
\[ \lambda(y^{-1})=  
 \left(
\begin{array}{cccccc}
 5 & -5 & 3 & -1 & 0 & 0 \\
 35 & -31 & 18 & -6 & 0 & 0 \\
 105 & -84 & 48 & -17 & 1 & 0 \\
 175 & -126 & 70 & -26 & 3 & 0 \\
 161 & -105 & 56 & -21 & 3 & 0 \\
 70 & -42 & 21 & -7 & 0 & 1 \\
\end{array}
\right) .
\]

The $\Lambda$-action on $\mathcal{T}_6$ is given by \[
t\cdot x = t\cdot \lambda(x^{-1}),
\] and
\[
t \cdot y =t \cdot \lambda(y^{-1}).
\]

Let $v_1 = (20, \ -10, \ 4, \ -1, \ 0, \ 0) \in \mathcal{T}_6$. One can easily check that $v_1\cdot y = -v_1 \cdot (1+x)$. Now, $[v_1)$ is a right $\Lambda$-module with $\Lambda$-action given by the representation $\lambda$. Also, since $v_1 \in R(1), [v_1)\subset R(1)$, but $R(1)$ is generated by $(1 ,\ 0, \ 0, \ 0, \ 0, \ 0 )$ and 
 
 \[  (20, \ -10, \ 4, \ -1, \ 0, \ 0)
  \left( {\begin{array}{cccccc}
  -1 & 0 & 0 & 0 & 0 & 0   \\
 0 & 0 & 0 & 0 & 0 & 0    \\ 
 0 & 0 & 0 & 0 & 0 & 0    \\
 0 & 0 & 0 & 0 & 0 & 0   \\
 -21 & 0 & 0 & 0 & 0 & 0   \\
  0  & 0 & 0 & 0 & 0 & 0  \\
  \end{array} } \right)=(1 ,\ 0, \ 0, \ 0, \ 0, \ 0 ),
\] therefore $[v_1)=R(1)$. We deduce that a six-dimensional $\Lambda$-lattice $M$ is isomorphic to $R(1)$ if and only if the following conditions are satisfied:
\begin{itemize}
\item $M(\Sigma)$: $m\cdot (1+x+x^2+\dots + x^6)=0$ for each $m \in M$;
\item $M(1)$: $M$ has a generator $v_1$ such that $v_1 \cdot y = -v_1 \cdot (1+x)$.
\end{itemize}

Similarly, if we define $v_2 = (7 ,\ -1, \ 0, \ 0, \ 0, \ 0 ) \in \mathcal{T}_6$, then it is easily checked that $v_2 \cdot y = v_2 \cdot (1+x)^2$ and $[v_2)=R(2)$.  We deduce that a six-dimensional $\Lambda$-lattice $M$ is isomorphic to $R(2)$ if and only if the following conditions are satisfied:
\begin{itemize}
\item $M(\Sigma)$: $m\cdot (1+x+x^2+\dots + x^6)=0$ for each $m \in M$;
\item $M(2)$: $M$ has a generator $v_2$ such that $v_2 \cdot y = v_2 \cdot (1+x)^2$.
\end{itemize}
Similarly, if we define $v_3 = (21, \ -7, \ 1, \ 1, \ -1, \ 0 )$, then $v_3 \cdot y = -v_3$ and $[v_3)=R(3)$. We deduce that a six-dimensional $\Lambda$-lattice $M$ is isomorphic to $R(3)$ if and only if the following conditions are satisfied:
\begin{itemize}
\item $M(\Sigma)$: $m\cdot (1+x+x^2+\dots + x^6)=0$ for each $m \in M$;
\item $M(3)$: $M$ has a generator $v_3$ such that $v_3 \cdot y = -v_3$.
\end{itemize}
Similarly, if we define $v_4 = (0, \ 0, \ 0, \ 1, \ -2, \ 2 )$, then $v_4 \cdot y = v_4\cdot(1+x)$ and $[v_4)=R(4)$.  We deduce that a six-dimensional $\Lambda$-lattice $M$ is isomorphic to $R(4)$ if and only if the following conditions are satisfied:
\begin{itemize}
\item $M(\Sigma)$: $m\cdot (1+x+x^2+\dots + x^6)=0$ for each $m \in M$;
\item $M(4)$: $M$ has a generator $v_4$ such that $v_4 \cdot y = v_4\cdot(1+x)$.
\end{itemize}
Similarly, if we define $v_5 = (77, \ -49, \ 28, \ -14, \ 6, \ -2 )$, then one can easily check that $v_5 \cdot y = -v_5\cdot(1+x)^2$ and $[v_5)=R(5)$.  We deduce that a six-dimensional $\Lambda$-lattice $M$ is isomorphic to $R(5)$ if and only if the following conditions are satisfied:
\begin{itemize}
\item $M(\Sigma)$: $m\cdot (1+x+x^2+\dots + x^6)=0$ for each $m \in M$;
\item $M(5)$: $M$ has a generator $v_5$ such that $v_5 \cdot y = -v_5\cdot(1+x)^2$.
\end{itemize}
Finally, if we define $v_6 = (7, \ -7, \ 7, \ -7, \ 7, \ -6 )$, then $v_6 \cdot y = v_6$ and $[v_6)=R(6)$.  We deduce that a six-dimensional $\Lambda$-lattice $M$ is isomorphic to $R(6)$ if and only if the following conditions are satisfied:
\begin{itemize}
\item $M(\Sigma)$: $m\cdot (1+x+x^2+\dots + x^6)=0$ for each $m \in M$;
\item $M(6)$: $M$ has a generator $v_6$ such that $v_6 \cdot y = v_6$.
\end{itemize}

Using the above characteristic equations, we can find representations \[\theta _k:\mathbb{Z}[G(7,6)]\rightarrow \mathcal{T}_6\] for each row $R(k)$ of $\mathcal{T}_6$. We give these explicitly now:

\[\theta_1(x^{-1})=\left(
\begin{array}{cccccc}
 0 & 0 & 0 & 0 & 0 & -1 \\
 1 & 0 & 0 & 0 & 0 & -1 \\
 0 & 1 & 0 & 0 & 0 & -1 \\
 0 & 0 & 1 & 0 & 0 & -1 \\
 0 & 0 & 0 & 1 & 0 & -1 \\
 0 & 0 & 0 & 0 & 1 & -1 \\
\end{array}
\right),\; \theta_1(y^{-1})= \left(
\begin{array}{cccccc}
 -1 & 1 & 0 & 0 & 0 & 0 \\
 -1 & 1 & 0 & -1 & 1 & 0 \\
 0 & 1 & 0 & -1 & 1 & 0 \\
 0 & 1 & -1 & 0 & 1 & 0 \\
 0 & 1 & -1 & 0 & 1 & -1 \\
 0 & 0 & 0 & 0 & 1 & -1 \\
\end{array}
\right), \]

\[\theta_2(x^{-1})=\left(
\begin{array}{cccccc}
 0 & 0 & 0 & 0 & 0 & -1 \\
 1 & 0 & 0 & 0 & 0 & -1 \\
 0 & 1 & 0 & 0 & 0 & -1 \\
 0 & 0 & 1 & 0 & 0 & -1 \\
 0 & 0 & 0 & 1 & 0 & -1 \\
 0 & 0 & 0 & 0 & 1 & -1 \\
\end{array}
\right),\; \theta_2(y^{-1})= \left(
\begin{array}{cccccc}
 1 & -1 & 0 & 0 & 1 & -1 \\
 2 & -2 & 0 & 1 & 0 & -1 \\
 1 & -2 & 0 & 2 & -1 & -1 \\
 0 & -2 & 1 & 1 & -1 & -1 \\
 0 & -2 & 2 & 0 & -1 & 0 \\
 0 & -1 & 1 & 0 & -1 & 1 \\
\end{array}
\right), \]

\[\theta_3(x^{-1})=\left(
\begin{array}{cccccc}
 0 & 0 & 0 & 0 & 0 & -1 \\
 1 & 0 & 0 & 0 & 0 & -1 \\
 0 & 1 & 0 & 0 & 0 & -1 \\
 0 & 0 & 1 & 0 & 0 & -1 \\
 0 & 0 & 0 & 1 & 0 & -1 \\
 0 & 0 & 0 & 0 & 1 & -1 \\
\end{array}
\right),\; \theta_3(y^{-1})= \left(
\begin{array}{cccccc}
 -1 & 0 & 0 & 0 & 1 & 0 \\
 0 & 0 & 0 & -1 & 1 & 0 \\
 0 & 0 & 0 & 0 & 1 & 0 \\
 0 & 0 & -1 & 0 & 1 & 0 \\
 0 & 0 & 0 & 0 & 1 & -1 \\
 0 & -1 & 0 & 0 & 1 & 0 \\
\end{array}
\right), \]

\[\theta_4(x^{-1})=\left(
\begin{array}{cccccc}
 0 & 0 & 0 & 0 & 0 & -1 \\
 1 & 0 & 0 & 0 & 0 & -1 \\
 0 & 1 & 0 & 0 & 0 & -1 \\
 0 & 0 & 1 & 0 & 0 & -1 \\
 0 & 0 & 0 & 1 & 0 & -1 \\
 0 & 0 & 0 & 0 & 1 & -1 \\
\end{array}
\right),\; \theta_4(y^{-1})= \left(
\begin{array}{cccccc}
 1 & -1 & 0 & 0 & 0 & 0 \\
 1 & -1 & 0 & 1 & -1 & 0 \\
 0 & -1 & 0 & 1 & -1 & 0 \\
 0 & -1 & 1 & 0 & -1 & 0 \\
 0 & -1 & 1 & 0 & -1 & 1 \\
 0 & 0 & 0 & 0 & -1 & 1 \\
\end{array}
\right), \]

\[\theta_5(x^{-1})=\left(
\begin{array}{cccccc}
 0 & 0 & 0 & 0 & 0 & -1 \\
 1 & 0 & 0 & 0 & 0 & -1 \\
 0 & 1 & 0 & 0 & 0 & -1 \\
 0 & 0 & 1 & 0 & 0 & -1 \\
 0 & 0 & 0 & 1 & 0 & -1 \\
 0 & 0 & 0 & 0 & 1 & -1 \\
\end{array}
\right),\; \theta_5(y^{-1})= \left(
\begin{array}{cccccc}
 -1 & 1 & 0 & 0 & -1 & 1 \\
 -2 & 2 & 0 & -1 & 0 & 1 \\
 -1 & 2 & 0 & -2 & 1 & 1 \\
 0 & 2 & -1 & -1 & 1 & 1 \\
 0 & 2 & -2 & 0 & 1 & 0 \\
 0 & 1 & -1 & 0 & 1 & -1 \\
\end{array}
\right), \]

\[\theta_6(x^{-1})=\left(
\begin{array}{cccccc}
 0 & 0 & 0 & 0 & 0 & -1 \\
 1 & 0 & 0 & 0 & 0 & -1 \\
 0 & 1 & 0 & 0 & 0 & -1 \\
 0 & 0 & 1 & 0 & 0 & -1 \\
 0 & 0 & 0 & 1 & 0 & -1 \\
 0 & 0 & 0 & 0 & 1 & -1 \\
\end{array}
\right),\; \theta_6(y^{-1})= \left(
\begin{array}{cccccc}
 1 & 0 & 0 & 0 & -1 & 0 \\
 0 & 0 & 0 & 1 & -1 & 0 \\
 0 & 0 & 0 & 0 & -1 & 0 \\
 0 & 0 & 1 & 0 & -1 & 0 \\
 0 & 0 & 0 & 0 & -1 & 1 \\
 0 & 1 & 0 & 0 & -1 & 0 \\
\end{array}
\right). \]
\end{subsection}

\begin{subsection}{The mapping $\pi_*:\Lambda \rightarrow [\pi)$}
In this section, we define an element $\pi \in \Lambda$ and show that $[\pi) \cong R(1)$ using the characteristic equations from section \ref{T_6}. We then go on to find the kernel of the map $\pi_*:\Lambda \rightarrow [\pi)$ defined by $\alpha \mapsto \pi \cdot \alpha$.

We begin by defining $$\pi = (x-1)((2+x^2+x^5)y+(-1+x^2+2x^3+2x^4+x^5)y^2+y^3)(1-y^3).$$

Let $v=\pi(1+x^2)$, Clearly $[v)\subset [\pi)$, but
\begin{align*} 
v(1+x)(1+x^4)&= \pi(1+x^2)(1+x)(1+x^4)\\
&= \pi(1+x+x^2 + \dots + x^7)\\
& = \pi
\end{align*}
and so $[v)=[\pi)$.
 One can check easily that $v\cdot y =-v\cdot(1+x)$, and so since $[\pi)$ is clearly a 6-dimensional $\Lambda$-lattice, we have shown:

\begin{proposition*}
The right $\Lambda$-module $[\pi)$ is isomorphic to $R(1)$.
\end{proposition*}

We have therefore constructed an explicit description of a surjection $\pi_* :  \Lambda \rightarrow R(1)$. Let $K = Ker(\pi_*)$. We now proceed to find a $\mathbb{Z}$-basis for $K$ which includes a $\mathbb{Z}$-basis for $R(1,3,4,5,6)$. This will be followed by a proof that $K/R(1,3,4,5,6) \cong \mathbb{Z}[C_6]$, and so $K$ is described by an extension of the type mentioned in section \ref{M(p)}. Utilising the identity $(1-y^3)(1+y^3)=0$, we see that $[y^3+1)\cap [x-1) \subset K$. We now define a $\mathbb{Z}$-basis $\{e(i)\}_{1\leq i \leq 18}$ for $[y^3+1)\cap [x-1)$:
\begin{itemize}
\item $e(i) = (y^3+1)(x-1)x^{i-1}$,  $1 \leq i \leq 6$;
\item $e(i) = (y^4+y)(x-1)x^{i-7}$,  $7 \leq i \leq 12$;
\item $e(i) = (y^5+y^2)(x-1)x^{i-13}$,  $13 \leq i \leq 18$.
\end{itemize}
If we define the representation $L^{\prime}:\Lambda \rightarrow GL_6(\mathbb{Z})$ by \[L^{\prime}(x^{-1})=\left(
\begin{array}{cccccc}
 0 & 0 & 0 & 0 & 0 & -1 \\
 1 & 0 & 0 & 0 & 0 & -1 \\
 0 & 1 & 0 & 0 & 0 & -1 \\
 0 & 0 & 1 & 0 & 0 & -1 \\
 0 & 0 & 0 & 1 & 0 & -1 \\
 0 & 0 & 0 & 0 & 1 & -1 \\
\end{array}
\right); \; L^{\prime}(y^{-1})=\left(
\begin{array}{cccccc}
 1 & 0 & 0 & 0 & 0 & 0 \\
 1 & 0 & -1 & 1 & 0 & 0 \\
 1 & 0 & -1 & 1 & 0 & -1 \\
 1 & -1 & 0 & 1 & 0 & -1 \\
 1 & -1 & 0 & 1 & -1 & 0 \\
 0 & 0 & 0 & 1 & -1 & 0 \\
\end{array}
\right),\]
standard calculations on the basis $\{e(i)\}_{1\leq i \leq 18}$ produce the representation $L:\Lambda \rightarrow GL_{18}(\mathbb{Z})$ for $[y^3+1)\cap [x-1)$ given by

\[L(x^{-1})=\left(
\begin{array}{ccc}
 L^{\prime}(x^{-1}) & 0 & 0 \\
 0 & L^{\prime}(x^{-1}) & 0 \\
 0 & 0 & L^{\prime}(x^{-1}) \\
\end{array}
\right);\] and 
\[ L(y^{-1})= \left(
\begin{array}{ccc}
 0 & 0 & L^{\prime}(y^{-1}) \\
 L^{\prime}(y^{-1}) & 0 & 0 \\
 0 & L^{\prime}(y^{-1}) & 0 \\
\end{array}
\right).\]
In order to show that $[y^3+1)\cap [x-1) \cong R(1,3,5)$, we begin by defining the following representation $\theta_{1,3,5}:\Lambda \rightarrow GL_{18}(\mathbb{Z})$ for $R(1,3,5)$ using the representations found in section \ref{T_6}:

\[\theta_{1,3,5}(x^{-1})=\left(
\begin{array}{ccc}
 \theta_1(x^{-1}) & 0 & 0 \\
 0 & \theta_3(x^{-1}) & 0 \\
 0 & 0 & \theta_5(x^{-1}) \\
\end{array}
\right);\] and 
\[\theta_{1,3,5}(y^{-1})=\left(
\begin{array}{ccc}
 \theta_1(y^{-1}) & 0 & 0 \\
 0 & \theta_3(y^{-1}) & 0 \\
 0 & 0 & \theta_5(y^{-1}) \\
\end{array}
\right).\] Now, if we define  \small\[ h =\left(
\begin{array}{cccccccccccccccccc}
 0 & 0 & -1 & 1 & 0 & 0 & 0 & 0 & 0 & 1 & -1 & 0 & 0 & 1 & -1 & 0 & 0 & 0 \\
 0 & 0 & -1 & 0 & 1 & 0 & 0 & 0 & 0 & 1 & 0 & -1 & 0 & 1 & 0 & -1 & 0 & 0 \\
 0 & 0 & -1 & 0 & 0 & 1 & 0 & 0 & 0 & 1 & 0 & 0 & 0 & 1 & 0 & 0 & -1 & 0 \\
 0 & 0 & -1 & 0 & 0 & 0 & -1 & 0 & 0 & 1 & 0 & 0 & 0 & 1 & 0 & 0 & 0 & -1 \\
 1 & 0 & -1 & 0 & 0 & 0 & 0 & -1 & 0 & 1 & 0 & 0 & 0 & 1 & 0 & 0 & 0 & 0 \\
 0 & 1 & -1 & 0 & 0 & 0 & 0 & 0 & -1 & 1 & 0 & 0 & -1 & 1 & 0 & 0 & 0 & 0 \\
 0 & 0 & -1 & 1 & 0 & 0 & 0 & -1 & 0 & 0 & 0 & 0 & 1 & -1 & 1 & -1 & 1 & 0 \\
 0 & 0 & -1 & 0 & 1 & 0 & 1 & -1 & -1 & 0 & 0 & 0 & 0 & 0 & 0 & 0 & 0 & 1 \\
 0 & 0 & -1 & 0 & 0 & 1 & 1 & 0 & -1 & -1 & 0 & 0 & 0 & -1 & 1 & -1 & 1 & 0 \\
 0 & 0 & -1 & 0 & 0 & 0 & 1 & 0 & 0 & -1 & -1 & 0 & 1 & -1 & 0 & 0 & 0 & 1 \\
 1 & 0 & -1 & 0 & 0 & 0 & 1 & 0 & 0 & 0 & -1 & -1 & 0 & 0 & 0 & -1 & 1 & 0 \\
 0 & 1 & -1 & 0 & 0 & 0 & 1 & 0 & 0 & 0 & 0 & -1 & 1 & -1 & 1 & -1 & 0 & 1 \\
 0 & 0 & -1 & 1 & 0 & 0 & 0 & 0 & -1 & 0 & 0 & 1 & -1 & 0 & 0 & 1 & -1 & 0 \\
 0 & 0 & -1 & 0 & 1 & 0 & 0 & 0 & -1 & -1 & 0 & 1 & 0 & -1 & 0 & 1 & 0 & -1 \\
 0 & 0 & -1 & 0 & 0 & 1 & 1 & 0 & -1 & -1 & -1 & 1 & 0 & 0 & -1 & 1 & 0 & 0 \\
 0 & 0 & -1 & 0 & 0 & 0 & 1 & 1 & -1 & -1 & -1 & 0 & -1 & 0 & 0 & 0 & 0 & 0 \\
 1 & 0 & -1 & 0 & 0 & 0 & 1 & 1 & 0 & -1 & -1 & 0 & 0 & -1 & 0 & 1 & -1 & 0 \\
 0 & 1 & -1 & 0 & 0 & 0 & 0 & 1 & 0 & 0 & -1 & 0 & 0 & 0 & -1 & 1 & 0 & -1 \\
\end{array}
\right)\]
Then $h$ has determinant $1$ and satisfies the following: 
\begin{itemize}
\item $h\cdot \theta_{1,3,5}(x^{-1})\cdot h^{-1} = L(x^{-1});$
\item $h\cdot \theta_{1,3,5}(y^{-1})\cdot h^{-1} = L(y^{-1}).$
\end{itemize}
We deduce that $[y^3+1)\cap [x-1) = [\eta(1)) \dotplus [\eta(3)) \dotplus [\eta(5))$ where $[\eta(i)) \cong R(i)$ and the elements $\eta(i)$ are defined as follows:
\begin{itemize}
\item $\eta(1) = h(v_1)=(y^3+1)(1+y+y^2)(x^5-x^4);$
\item $\eta(3) = h(v_3)=(y^3+1)(-(x^4-x^3)+y(x^6-x)+y^2(x^5-x^2));$
\item $\eta(5) =h(v_5)= (y^3+1)(-(x^6-x^5)+y((x^6-x^5)+(x^4-x^3)+(x-1))+y^2(-(x^4-x^3)-(x-1)).$
\end{itemize}
Now, as $\{e(i)\}_{1\leq i \leq 18} \cup \{y^i \cdot x^j \mid 0\leq i \leq 2,\; 0 \leq j \leq 6\}\cup \{ y^3,y^4,y^5 \}$ is a basis for $\Lambda$, and $\{ e(i) \}_{1 \leq i \leq 18}$ is a basis for $[y^3+1)\cap [x-1)$, $\Lambda/[y^3+1)\cap [x-1)$ is torsion free. From the exact sequence
 $$0 \rightarrow K/[y^3+1)\cap [x-1) \rightarrow \Lambda/[y^3+1)\cap [x-1) \rightarrow R(1) \rightarrow 0,$$ we deduce that $K/[y^3+1)\cap[x-1)$ is torsion free. We therefore have the following short exact sequence of $\Lambda$-lattices: $$0\rightarrow [y^3+1)\cap [x-1) \rightarrow K\rightarrow K/[y^3+1)\cap [x-1) \rightarrow 0.$$ We can now form a basis for $K$ using bases for $[y^3+1)\cap [x-1)$ and $K/[y^3+1)\cap [x-1)$. As $[y^3+1)\cap [x-1)\cong R(1,3,5)$, in order to find a basis for $K$ which contains a basis for $R(1,3,4,5,6)$, we will find elements $\eta(4),\eta(6)\in K$ such that there exists a basis for $K/[y^3+1)\cap[x-1)$ which contains the set $\{\eta(i)X^j\mid i=4,6 ;\;  j=0,1,\dots ,5 \}$ and  $[\eta(i))\cong R(i)$ for $i=4,6$. Here, we use capitalisation to represent the image of $x$ in $K/[y^3+1)\cap[x-1)$. We define $\eta(4), \eta(6)$ as follows,
\begin{itemize}
\item $\eta(4)=(x-1)((1+x^5)+(-x+x^4+x^5)y+(x^2+x^3+x^4+x^5)y^2+(x^3+x^5)y^3+(-1-x-2x^2-x^3-x^4)y^4+(-x^2-x^3-x^4)y^5)$,
\item $\eta(6)=(x-1)(1+x^5)(1+y+y^2+y^3+y^4+y^5)$.
\end{itemize}
Through tedious calculations one can check that $\pi\cdot \eta(4)=0$ and $\pi \cdot \eta(6)=0$ and so $\eta(4),\eta(6)\in K$; it is immediately clear that $\eta(6)\cdot y=\eta(6)$, and through another tedious calculation one can check that $\eta(4)\cdot y=\eta(4)\cdot (1+x)$, therefore $[\eta(i))\cong R(i)$ for $i=4,6$. We will now find a basis for $K/[y^3+1)\cap[x-1)$, which we transform to a basis containing $\{\eta(i)X^j \mid i=4,6 \; j=0,1,\dots ,5\}$. We begin by expressing the mapping $$\overline{\pi_*}:\Lambda/[y^3+1)\cap[x-1) \rightarrow [\pi), \;\;\; \; \alpha + [y^3+1)\cap [x-1) \mapsto \pi \cdot \alpha$$ in matrix form, where we take $\Lambda/[y^3+1)\cap[x-1)$ to have basis 
\begin{align*}
 \{\,&1,X,X^2,X^3,X^4,X^5,X^6,Y,YX,YX^2,YX^3,YX^4,YX^5,YX^6 \\
  & Y^2,Y^2X,Y^2X^2,Y^2X^3,Y^2X^4,Y^2X^5,Y^2X^6,Y^3,Y^4,Y^5 \}.
\end{align*}
Here, capitalisation is used to represent the image in $\Lambda/[y^3+1)\cap[x-1)$. We take $[\pi)$ to have basis \[ \{\pi,\pi\cdot x ,\pi\cdot x^2 ,\pi\cdot x^3 ,\pi\cdot x^4 ,\pi\cdot x^5\}.\] Now,  $\pi(1+x^2)y=-\pi(1+x^2)(1+x)$ and so 
\begin{align*}
\pi y &=\pi(1+x^2)(1+x)(1+x^4) y\\ &=\pi(1+x^2)y(1+x^5)(1+x^6)\\ &=-\pi\cdot(1+x^2)(1+x)(1+x^5)(1+x^6),\\
&=\pi(-1+x^2+2x^3+2x^4+x^5)
\end{align*} Using this equality, we can form the following matrix for $\overline{\pi_*}$:
\small\[\arraycolsep=2pt\def\arraystretch{1.2}\left(
\begin{array}{cccccccccccccccccccccccc}
 1 & 0 & 0 & 0 & 0 & 0 & -1 & -1 & -1 & -1 & 0 & 1 & 1 & 1 & 2 & -1 & 1 & 0 & -1 & 1 & -2 & -1 & 1 & -2 \\
 0 & 1 & 0 & 0 & 0 & 0 & -1 & 0 & -2 & -2 & -1 & 1 & 2 & 2 & 0 & 1 & 0 & 1 & -1 & 0 & -1 & 0 & 0 & 0 \\
 0 & 0 & 1 & 0 & 0 & 0 & -1 & 1 & -1 & -3 & -2 & 0 & 2 & 3 & 1 & -1 & 2 & 0 & 0 & 0 & -2 & 0 & -1 & -1 \\
 0 & 0 & 0 & 1 & 0 & 0 & -1 & 2 & 0 & -2 & -3 & -1 & 1 & 3 & 0 & 0 & 0 & 2 & -1 & 1 & -2 & 0 & -2 & 0 \\
 0 & 0 & 0 & 0 & 1 & 0 & -1 & 2 & 1 & -1 & -2 & -2 & 0 & 2 & 0 & -1 & 1 & 0 & 1 & 0 & -1 & 0 & -2 & 0 \\
 0 & 0 & 0 & 0 & 0 & 1 & -1 & 1 & 1 & 0 & -1 & -1 & -1 & 1 & 1 & -1 & 0 & 1 & -1 & 2 & -2 & 0 & -1 & -1 \\
\end{array}
\right).\]

By utilising elementary linear algebra we can now easily find a $\mathbb{Z}$-basis for $K/[y^3+1)\cap[x-1)$. By expressing $\{\eta(i)X^j\mid i=4,6 \;, j=0,1,\dots ,5 \}$ in terms of the basis for Ker$(\overline{\pi_*})$ and utilising the Smith Normal Form, we construct the following basis for $K/[y^3+1)\cap[x-1)$:
 \begin{align*}\{\,& \eta(4)X^i \mid 0\leq i \leq 6 \} \cup \{ \eta(6)X^i \mid 0\leq i \leq 6 \} \cup\{1+Y^3,-2-X^2-X^5+Y^2,\\ &2+X^2+X^5+Y^5,-1+X^2+2X^3+2X^4+X^5+Y^4, 
   1-X^2-2X^3-2X^4-X^5+Y,\\&1+X+X^2+X^3+X^4+X^5+X^5 \}.
\end{align*}
We therefore have a basis for $K$ which includes bases for $[\eta(i))$ for $i=1,3,4,5,6$. One can now calculate the representation $\rho: \Lambda \rightarrow GL_6(\mathbb{Z})$ for $K/R(1,3,4,5,6)$.
\[\rho(x^{-1})=\left(
\begin{array}{cccccc}
 1 & 0 & 0 & 0 & 0 & 0 \\
 0 & 1 & 0 & 0 & 0 & 0 \\
 0 & 0 & 1 & 0 & 0 & 0 \\
 0 & 0 & 0 & 1 & 0 & 0 \\
 0 & 0 & 0 & 0 & 1 & 0 \\
 0 & 0 & 0 & 0 & 0 & 1 \\
\end{array}
\right); \; \rho(y^{-1})= \left(
\begin{array}{cccccc}
 0 & 1 & 0 & 0 & 0 & 0 \\
 0 & 0 & 0 & 0 & 1 & 0 \\
 0 & 0 & 0 & 1 & 0 & 0 \\
 1 & 0 & 0 & 0 & 0 & 0 \\
 1 & -4 & 4 & 5 & -5 & 7 \\
 0 & -3 & 3 & 3 & -3 & 5 \\
\end{array}
\right). \]
Now, if we let \[f=\left(
\begin{array}{cccccc}
 0 & 1 & 1 & 0 & 0 & -1 \\
 1 & 1 & 0 & 0 & -1 & 0 \\
 0 & -1 & 0 & 1 & 1 & 0 \\
 -1 & 0 & 1 & 1 & 0 & 0 \\
 1 & 0 & 0 & -1 & 0 & 1 \\
 2 & 1 & -1 & -2 & -1 & 1 \\
\end{array}
\right).\] Then $f^{-1}\rho(g)f=\sigma(g)$ for each $g\in \Lambda$, where $\sigma(g)$ is the regular representation of the $\Lambda$-module $\mathbb{Z}[C_6]$. We deduce that $K$ lies in a short exact sequence of the form $$0\rightarrow R(1,3,4,5,6)\rightarrow K \rightarrow \mathbb{Z}[C_6]\rightarrow 0.$$ By \textbf{Theorem \ref{k-invariants}}, to prove that $K\cong \Lambda/i(R(2))$ for some injective $\Lambda$-homomorphism $i:R(2)\rightarrow \Lambda$, it is sufficient to show that the above exact sequence has all non-zero $k$-invariants. The representation $\varphi:\Lambda \rightarrow GL_{36}(\mathbb{Z})$ of $K$ given by the above exact sequence after the change of basis defined by $f$ takes the following form:
\[\varphi(x^{-1})=\left(
\begin{array}{cccccc}
 \theta_1(x^{-1}) & 0 & 0 & 0 & 0 & C(1)  \\
 0 & \theta_3(x^{-1}) & 0 & 0 & 0 & C(3) \\
 0 & 0 & \theta_4(x^{-1}) & 0 & 0 & C(4) \\
 0 & 0 & 0 & \theta_5(x^{-1}) & 0 & C(5)\\
 0 & 0 & 0 & 0 & \theta_6(x^{-1}) & C(6) \\
 0 & 0 & 0 & 0 & 0 & I_6 \\

\end{array}
\right);\]

\[\varphi(y^{-1})=\left(
\begin{array}{cccccc}
 \theta_1(y^{-1}) & 0 & 0 & 0 & 0 & D(1)  \\
 0 & \theta_3(y^{-1}) & 0 & 0 & 0 & D(3) \\
 0 & 0 & \theta_4(y^{-1}) & 0 & 0 & D(4) \\
 0 & 0 & 0 & \theta_5(y^{-1}) & 0 & D(5)\\
 0 & 0 & 0 & 0 & \theta_6(y^{-1}) & D(6) \\
 0 & 0 & 0 & 0 & 0 & \sigma(y^{-1})\\
\end{array}
\right).\]

Now, to check that the $k$-invariant corresponding to $Ext^1_{\Lambda}(\mathbb{Z}[C_{6}],R(i))$ for $i=1,3,4,5,6$ is non-zero, we must show that the extension defined by 
\[\varphi_i(x^{-1})=\left(
\begin{array}{cc}
 \theta_i(x^{-1}) & C(i) \\
 0 & I_6 \\
\end{array}
\right); \;
\varphi_i(y^{-1})=\left(
\begin{array}{cc}
 \theta_i(y^{-1}) & D(i) \\
 0 & \sigma(y^{-1}) \\
\end{array}
\right),\] is not congruent to the trivial extension in $Ext^{1}_{\Lambda}(R(i),\mathbb{Z}[C_6])$ for any $i$. This is equivalent to showing that there is no matrix $\psi_i \in GL_{12}(\mathbb{Z})$ of the form \[\psi_i=\left(
\begin{array}{cc}
 I_6 & X_i \\
 0 & I_6 \\
\end{array}
\right)\] such that
\[\left(
\begin{array}{cc}
 I_6 & X_i \\
 0 & I_6 \\
\end{array}
\right)\left(
\begin{array}{cc}
 \theta_i(x^{-1}) & C(i) \\
 0 & I_6 \\
\end{array}
\right)=\left(
\begin{array}{cc}
 \theta_i(x^{-1}) & 0 \\
 0 & I_6 \\
\end{array}
\right)\left(
\begin{array}{cc}
 I_6 & X_i \\
 0 & I_6 \\
\end{array}
\right),\]
and

\[\left(
\begin{array}{cc}
 I_6 & X_i \\
 0 & I_6 \\
\end{array}
\right)\left(
\begin{array}{cc}
 \theta_i(y^{-1}) & D(i) \\
 0 & \sigma(y^{-1}) \\
\end{array}
\right)=\left(
\begin{array}{cc}
 \theta_i(y^{-1}) & 0 \\
 0 & \sigma(y^{-1}) \\
\end{array}
\right)\left(
\begin{array}{cc}
 I_6 & X_i \\
 0 & I_6 \\
\end{array}
\right).\]
This is equivalent to showing that there is no $X_i\in M_{6\times 6}(\mathbb{Z})$ such that

\begin{equation}\label{important eq} C(i)=(\theta_i(x^{-1})-I_6)X_i
\end{equation} and
\begin{equation}
 D(i)+X_i\sigma(y^{-1})=\theta_i(y^{-1})X_i
\end{equation}
In our case, we calculate the $C(i)$ and $D(i)$ to be:
\[C(1)=\left(
\begin{array}{cccccc}
 -5 & -2 & 2 & 5 & 2 & -2 \\
 -4 & 1 & 2 & 4 & -2 & -1 \\
 -5 & -3 & 2 & 4 & 3 & -1 \\
 -1 & 0 & 1 & 1 & 0 & 0 \\
 0 & -1 & 0 & -1 & 1 & 1 \\
 0 & -3 & -1 & 0 & 2 & 2 \\
\end{array}
\right);\;
C(3)=\left(
\begin{array}{cccccc}
 7 & 3 & -2 & -8 & -1 & 1 \\
 5 & -1 & -2 & -6 & 3 & 1 \\
 14 & 4 & -6 & -14 & -4 & 6 \\
 7 & -5 & -4 & -8 & 7 & 3 \\
 8 & 4 & -3 & -7 & -3 & 1 \\
 3 & 3 & 0 & -4 & -1 & -1 \\
\end{array}
\right);\] \[C(4)=\left(
\begin{array}{cccccc}
 0 & 0 & 0 & 0 & 0 & 0 \\
 0 & 0 & 0 & 0 & 0 & 0 \\
 1 & -2 & -1 & -1 & 2 & 1 \\
 2 & 0 & -1 & -2 & 0 & 1 \\
 1 & -2 & -1 & -1 & 2 & 1 \\
 0 & 0 & 0 & 0 & 0 & 0 \\
\end{array}
\right);\;
C(5)=\left(
\begin{array}{cccccc}
 -7 & 0 & 4 & 6 & 1 & -4 \\
 -4 & 3 & 3 & 3 & -3 & -2 \\
 0 & -2 & 0 & -2 & 2 & 2 \\
 7 & -3 & -4 & -8 & 3 & 5 \\
 10 & -2 & -5 & -11 & 3 & 5 \\
 5 & -2 & -3 & -5 & 2 & 3 \\
\end{array}
\right);\]
\[C(6)=\left(
\begin{array}{cccccc}
 2 & 0 & -1 & -2 & 0 & 1 \\
 2 & 0 & -1 & -2 & 0 & 1 \\
 0 & 0 & 0 & 0 & 0 & 0 \\
 0 & 0 & 0 & 0 & 0 & 0 \\
 -1 & 2 & 1 & 1 & -2 & -1 \\
 0 & 0 & 0 & 0 & 0 & 0 \\
\end{array}
\right).\]

\[D(1)=\left(
\begin{array}{cccccc}
 1 & -1 & 2 & -1 & 1 & -2 \\
 1 & 0 & 5 & -1 & 0 & -5 \\
 2 & 2 & 6 & -2 & -2 & -6 \\
 2 & 1 & 5 & -2 & -1 & -5 \\
 0 & -1 & 3 & 0 & 1 & -3 \\
 1 & -1 & 0 & -1 & 1 & 0 \\
\end{array}
\right);\;
D(3)=\left(
\begin{array}{cccccc}
 -7 & 3 & 4 & 7 & -3 & -4 \\
 -2 & -2 & 0 & 2 & 2 & 0 \\
 -7 & -2 & -3 & 7 & 2 & 3 \\
 -6 & -1 & -7 & 6 & 1 & 7 \\
 -4 & -3 & -7 & 4 & 3 & 7 \\
 -9 & -2 & -1 & 9 & 2 & 1 \\
\end{array}
\right);\] \[D(4)=\left(
\begin{array}{cccccc}
 0 & 1 & 1 & 0 & -1 & -1 \\
 1 & 0 & 1 & -1 & 0 & -1 \\
 -1 & 0 & 1 & 1 & 0 & -1 \\
 1 & 0 & -1 & -1 & 0 & 1 \\
 -1 & 0 & -1 & 1 & 0 & 1 \\
 0 & -1 & -1 & 0 & 1 & 1 \\
\end{array}
\right);\;
D(5)=\left(
\begin{array}{cccccc}
 3 & 2 & 7 & -3 & -2 & -7 \\
 4 & 3 & 13 & -4 & -3 & -13 \\
 3 & 4 & 13 & -3 & -4 & -13 \\
 4 & 2 & 6 & -4 & -2 & -6 \\
 -1 & -1 & 0 & 1 & 1 & 0 \\
 1 & -3 & -4 & -1 & 3 & 4 \\
\end{array}
\right);\]
\[D(6)=\left(
\begin{array}{cccccc}
 -1 & 0 & -1 & 1 & 0 & 1 \\
 -2 & 0 & -2 & 2 & 0 & 2 \\
 0 & -1 & -3 & 0 & 1 & 3 \\
 -3 & 0 & -1 & 3 & 0 & 1 \\
 1 & 0 & -1 & -1 & 0 & 1 \\
 -2 & 1 & 1 & 2 & -1 & -1 \\
\end{array}
\right).\]

One can check easily that there are no solutions $X_i\in M_{6\times 6}(\mathbb{Z})$ for $i=1,2,\dots ,6$ for equation $(\ref{important eq})$. We have shown that the kernel $K$ of the map $\pi_*$, lies in an extension of the form
\[
0 \rightarrow R(1,3,4,5,6) \rightarrow K \rightarrow \mathbb{Z}[C_{p-1}] \rightarrow 0
\] with all non-zero $k$-invariants. Therefore, by \textbf{Theorem \ref{k-invariants}}, and our discussion at the end of section \ref{M(p)}, we have shown:

\begin{theorem}\label{surjection}
Over $\Lambda=\mathbb{Z}[G(7,6)]$, $\Omega_3(\mathbb{Z})=[R(2)\oplus [y-1)]$ i.e. the condition $M(7)$ holds.
\end{theorem}
Therefore, by \textbf{Theorem \ref{If M(p) holds}} 

\begin{theorem}\label{syzygy}
The D(2)-property holds for $G=G(7,6)$.
\end{theorem}
\end{subsection}
\end{section}
\section{Acknowledgements}
The author would like to thank his PhD supervisor, Professor FEA Johnson, who suggested the problem, and who through their knowledge and insight has assisted the author greatly in finding its eventual solution.

\end{document}